\documentclass[11pt]{amsart}
\usepackage{amssymb,latexsym,amsfonts,verbatim,amscd,ifthen}
\usepackage[normalem]{ulem}
\usepackage[usenames,dvipsnames]{pstricks}
\usepackage{pstricks}
\usepackage{epsfig}
\usepackage{pst-grad} 
\usepackage{pst-plot}
\usepackage{graphics,hyperref} 

\newtheorem{thm1}{Theorem}[section]
\newtheorem{lem1}[thm1]{Lemma}

\newtheorem{rem1}[thm1]{Remark}

\newtheorem{def1}[thm1]{Definition}

\newtheorem{cor1}[thm1]{Corollary}

\newtheorem{prop1}[thm1]{Proposition}
\newtheorem{ex1}[thm1]{Example}
\newtheorem{ex2}[thm1]{Examples}

\DeclareMathOperator{\supp}{supp}

\DeclareMathOperator{\rank}{rank}

\newcommand{\A}{\mathcal{A}}

\newcommand{\N}{\mathbb{N}}
\newcommand{\Z}{\mathbb{Z}}

\newcommand{\Q}{\mathbb{Q}}
\newcommand{\B}{\mathbf{B}}
                                                                                                                                                                   \newcommand{\ub}{\bold u}
\newcommand{\vb}{\bold v}
\newcommand{\wb}{\bold w}
\newcommand{\zb}{\bold z}
\newcommand{\cb}{\bold c}

\newcommand{\ab}{\bold a}
\newcommand{\bb}{\bold b}

\begin{document}

\title[Minimal generating sets of lattice Ideals]{Minimal generating sets of lattice Ideals}
\author[H. Charalambous]{Hara Charalambous}
\address{ Department of Mathematics, Aristotle University of Thessaloniki, Thessaloniki \\54124, Greece} \email{hara@math.auth.gr}
\author[A. Thoma]{ Apostolos Thoma}
\address{ Department of Mathematics, University of Ioannina, Ioannina 45110, Greece } \email{athoma@uoi.gr}
\author[M. Vladoiu]{Marius Vladoiu}
\address{ Faculty of Mathematics and Computer Science, University of Bucharest,  Str. Academiei 14, Bucharest,
RO-010014, Romania} \email{vladoiu@fmi.unibuc.ro}

\begin{abstract} 
Let $L\subset \Z^n$ be a lattice and $I_L=\langle x^{\ub}-x^{\vb}:\ {\ub}-{\vb}\in L\rangle $ be the corresponding lattice ideal in $\Bbbk[x_1,\ldots, x_n]$, where $\Bbbk$ is a field.  In this paper we describe minimal binomial generating sets of $I_L$ and their invariants. We use as a main tool a graph construction on equivalence classes of fibers of $I_L$. As one application of the theory developed we characterize binomial complete intersection lattice ideals, a longstanding open problem in the case of non-positive lattices.
\end{abstract}
\maketitle

\section{Introduction}
Let $R=\Bbbk[x_1,\ldots, x_n]$ where $\Bbbk$ is a field, and let $L$ be a lattice in $\mathbb{Z}^n$. The lattice ideal $I_L$ is the ideal:
\[ 
I_L:=\langle x^{\ub}-x^{\vb}:\ {\ub}-{\vb}\in L\rangle\; .
\]
Throughout this paper we call the set $S$ a {\it minimal generating set} of $I_L$ if $S$ consists of binomials and is a minimal \footnote {The notion of Markov basis was introduced in \cite{DS} to denote  a generating  binomial set of $I_L$, see \cite[Theorem 3.1]{DS}. Occasionally in the literature, especially in algebraic statistics, the minimal generating sets are also referred to as minimal Markov bases.} system of generators of $I_L$. If the  minimal generating set $S$ of $I_L$ is of minimal cardinality, then we say that $S$  is a {\it cardinality-minimal generating set} and we write  $\mu(I_L)$  for $|S|$. \footnote {Even though for positive lattices the  notions of  
 minimal generating sets and cardinality-minimal generating sets coincide, this is not true for  non-positive lattices.} Lattice ideals were first systematically studied in \cite{ES}. We note that toric ideals are lattice ideals $I_L$ for which the lattice $L$ is the kernel of an integer matrix.  The study of lattice ideals is a rich subject on its own, see \cite{MS, St} for the general theory and \cite{EM} for recent developments. Moreover lattice ideals have applications in various areas of mathematics, such as algebraic statistics \cite{DS, HO}, integer programming \cite{DSS}, hypergeometric differential equations \cite{DMM}, graph theory \cite{Hi-O}, etc.  

There are two types of lattices in $\Z^n$, the positive ones and the non-positive ones. We say that $L$ is a {\it positive} lattice if $L\cap \N^n=\{\bf{0}\}$. Otherwise, 
we say that $L$ is  a  {\it non-positive} lattice. 
Almost all results in the literature deal with positive lattices, with  few exceptions like in \cite{CTV, ES, GRV, HM,KM, LV, OPV}. We survey below several well known facts for positive lattices, 
as they pertain to  our study. For the proofs we refer to \cite{MS} or \cite{St} unless otherwise noted. 
Let $L$ be a positive lattice. We let $\{{\bf e}_i:\  1\leq i\leq n\}$ be  the canonical basis of 
$\mathbb{Z}^n$  and $\A$ be the subsemigroup of $\mathbb{Z}^n/L$ generated by the elements  $\{{\bf a}_i={\bf e}_i+L:1\leq i\leq n\}$. 
The semigroup $\A$ is partially ordered:
\[
{\bf c} \geq {\bf d} \Longleftrightarrow \ \textrm{there is} \
{\bf e} \in \A \ \textrm{such that} \ {\bf c}={\bf d}+{\bf e}\; .
\]
For $x^{\vb}=x_1^{v_1} \cdots x_n^{v_n}$ 
we set 
\[ 
\deg_{ \A}(x^{\vb}):=v_1{\bf a}_1+\cdots+v_n{\bf a}_n \in \mathcal{A}\; .  
\] 
It follows that $I_L$ is $\A$-graded and that 
\[ 
I_L=\langle x^{\ub}-x^{\vb}:\ \deg_{ \A}(x^{\ub})=\deg_{ \A}(x^{\vb})\ \rangle\;. 
\]
For  ${\ub}\in \N^n$,
the $I_L$-{\it fiber} $F_{\ub}$ is the set \[ F_{\ub}:=\{ x^{\vb}:\ x^{\vb}-x^{\ub}\in I_L\}=\{ x^{\vb}:\ \deg_{\A} (x^{\vb})=\deg_{\A} (x^{\ub})\}\;.\] 
$F_{\ub}$ is a finite set for all $\ub\in \N^n$.  The homogeneous Nakayama Lemma applies and thus all minimal generating sets of $I_L$  have the same cardinality and are 
 cardinality-minimal generating sets. Let $S$ be a minimal generating set of $I_L$ and consider the multiset of  the
   fibers $F_{\ub}$, where $x^{\ub}$ is a monomial term of a binomial in $S$. 
 By the results of \cite{CKT} this multiset  is an invariant of $I_L$. 
 Lastly we note that the union{\footnote{The union of all  minimal generating sets of $I_L$ is called the universal Markov basis of $I_L$, following \cite[Definition 3.1]{HS}.}}  of all minimal generating sets of $I_L$ is a finite set, since all minimal generating sets are subsets of a finite set, the Graver basis of $I_L$, see \cite{St}.

\par How do these facts translate to non-positive lattices? Suppose  that $L$ is a non-positive lattice and let $A$ be (as before) the semigroup of $\mathbb{Z}^n/L$ 
generated by the elements  $\{{\bf a}_i={\bf e}_i+L:1\leq i\leq n\}$.
The ideal $I_L$ is    generated by the binomials $x^{\ub}-x^{\vb}$, with  $\deg_{ \A}(x^{\ub})=\deg_{ \A}(x^{\vb})$  and the $I_L$-fibers are defined in the same way
as in the positive case. 
However, the 
 semigroup $\A$ is no longer partially ordered and the $I_L$-fibers are no longer finite. There are minimal generating sets which are not cardinality-minimal generating sets. In the  nonprincipal case, the multiset of   the $I_L$-fibers that correspond to a given cardinality-minimal generating set is not an invariant of $I_L$ and the union of all cardinality-minimal generating sets 
of $I_L$ is an infinite set. 
 
Consider, for example,  the non-positive lattice $L$ of $\Z^2$ generated by  $(1,1), (5,0)$.  The ideal $I_L$ is not  principal and the sets $\{1-xy, 1-x^5 \}$, 
$\{ 1-xy, x^3-y^2\}$ are examples of cardinality-minimal generating sets of $I_L$. Thus $\mu(I_L)=2$. However $\{1-x^2y^2, 1-x^3y^3, 1-x^5\}$ is also a minimal generating set of $I_L$ and in fact  $I_L$ has a minimal generating set{\footnote{Let $p_1,\ldots, p_s$ be $s$ distinct primes and let $a_i=(p_1\cdots p_s)/p_i$.  Then  $\langle 1-(xy)^{a_1}, \ldots, 1-(xy)^{a_s}\rangle =\langle 1-xy\rangle$ and  $ \{ 1-x^5, 1-(xy)^{a_1},\ldots, 1-(xy)^{a_s}\}$ is a minimal generating set of $I_L$ of cardinality $s+1$.}} of cardinality $k$, $\forall k\ge 2$.{\footnote{There are cardinality-minimal generating sets that are not contained in the Graver basis of $I_L$,  as is the case for $\{1-x^{2012}y^{2017}, y^4-x^{2013}y^{2022}\} $. The union of the cardinality-minimal generating sets of $I_L$  is an infinite set.}} The $I_L$-fibers that contain $1$ and $x^3$ are distinct and thus the multiset of the fibers for
 the above cardinality-minimal generating sets of $I_L$ are also distinct, see Example~\ref{example_multi}(b) for a precise description of the $I_L$-fibers.

The main purpose of this paper is to determine  the invariants of the minimal generating sets of $I_L$, for any lattice $L\subset \Z^n$, and  to show how to compute cardinality-minimal generating sets of $I_L$. To do so, we     consider an equivalence relation among the  $I_L$-fibers and  order  the equivalence classes. We  appropriately choose
the binomial generators that correspond to the smallest possible equivalence class. Then we apply a graph construction to the equivalence classes and  determine trees on these graphs in order to identify the binomials that could be part of a minimal generating set. In the process we identify the invariants of the cardinality-minimal generating sets of $I_L$. As an application of the theory developed we characterize    all binomial complete intersection lattice ideals, a longstanding open problem. In more detail, the structure and the main results of this paper are as follows. 

In Section 2 we introduce and study  the sublattice $L_{pure}$ of $L$. This  is the sublattice of $L$  generated by $L\cap\N^n$ and in essence it measures the deviation of $L$ from being positive. For example, the fibers of $I_L$ are finite if and only if $L_{pure}$ is the zero lattice. The support $\sigma$ of $L_{pure}$ is crucial to our study. Here $\sigma$ is a subset of $[n]=\{1,\ldots, n\}$ and $i\in \sigma$ if and only if there is a $\ub=(u_1,\ldots, u_n)\in L_{pure}$ so that $u_i\neq 0$. In Corollary~\ref{cor_basis_pure} we show that a basis of $L_{pure}$ can be chosen so that its elements are in $\N^n$ and have full support $\sigma$. This will be needed in Section 4 when computing a basis of $I_{L_{pure}}$.
 
In Section 3 we define an equivalence relation among the $I_L$-fibers. According to Definition~\ref{equiv}, two $I_L$-fibers $F$, $G$ are equivalent if there exist $\ub,\vb\in \N^n$ such that $x^{\ub}F\subset G$ and $x^{\vb} G\subset F$. It is essential in our study that the cardinality of the equivalence classes of the $I_L$-fibers is a constant determined by $\sigma$, as Proposition~\ref{same_nr_equiv} and Lemma~\ref{finite_equivalence_fiber} show. Next we introduce a partial order `` $\leq_{_{I_L}}$" on the equivalence classes. This order is compatible with the $\A$-order in the case of a positive lattice $L$, and the set of equivalence classes has a smallest element, namely the equivalence class of the fiber that contains $1$. In Theorem~\ref{chain_least} we show that any strictly descending chain of equivalence classes terminates, thus allowing inductive arguments to work. Among all $I_L$-fibers, we are primarily interested in those whose equivalence classes correspond to elements participating in a cardinality-minimal generating set of $I_L$. These are called Markov fibers, see Definition~\ref{def_markov_fiber}, and in Corollary~\ref{markov_fib_inv_cor} we show that the set of these equivalence classes is an invariant of $I_L$.  

Section 4 is the core of this paper. Its main objective is the characterization of all cardinality-minimal generating sets of $I_L$. First we attack the case when $L=L_{pure}$: note that  in this case we know that $I_L$ is a binomial complete intersection, see \cite[Theorem 2.1]{ES}.  In the first main result of this section,  Theorem~\ref{Markov_pure}, we describe  all cardinality-minimal generating sets of $I_L$. We remark that if $L$ is not cyclic, then there are infinitely many such bases. Next we make the transition to the general non-positive  lattices  by taking a closer look to the equivalence classes not containing the fiber of $1$. If $\overline{F}$ is such an equivalence class, Definition~\ref{gamma_graph} constructs a graph $\Gamma_{\overline F}$ : the edges of $\Gamma_{\overline F}$ determine binomials that belong to minimal generating sets of $I_L$. We reach the main result of this section in Theorem~\ref{main_thm} which describes all cardinality-minimal generating sets of $I_L$, for $L$  an arbitrary lattice. We note that when $L$ is positive, the description specializes to the results of \cite{CKT}, see Remark~\ref{positive_main_thm} for more details. As a consequence of our main theorem we can compute $\mu(I_L)$. This computation is the point of Theorem~\ref{cor_rank+graph}. Next we refine the  conclusion of Corollary~\ref{markov_fib_inv_cor} about the invariants of $I_L$. Let $S$ be a cardinality-minimal generating set of $I_L$. We show that the  multiset, consisting of the equivalence classes of the Markov fibers determined by $S$, is  an invariant of all cardinality-minimal generating sets of $I_L$, see 
Corollary~\ref{multiset_inv}. This is the best result  one can expect, as the example of  $L=\langle (1,1), (0,5)\rangle$ points out. Finally, in  Corollary~\ref{cor_multiset} we give an invariant of all minimal generating sets of $I_L$. They all determine the same multiset of  equivalence classes of fibers away from the fiber of $1$. In the last part of Section 4 we discuss the indispensable binomials and monomials of $I_L$ and show that in the general case of non-positive lattices, the union of the cardinality-minimal generating sets of $I_L$ is an infinite set.  

In Section 5 we determine the lattices $L$ for which $I_L$ is a  binomial complete intersection  ideal,  see the beginning of the section for a short history pertaining to the case of positive lattices. We note that this problem  was completely solved for positive lattices, see \cite[Theorem 3.9]{MT}. We refer to \cite{MT} for a more comprehensive description of the problem.  Even though for positive lattices, the class of complete intersection ideals is the same as the class of complete intersection binomial ideals, this is an open problem for non-positive lattices. In Theorem~\ref{compint} we show that $I_L$ is a binomial complete intersection ideal if and only if $I_{L^{\sigma}}$ is a complete intersection ideal, where $L^\sigma$ is the lattice in the complement of $\sigma$. In Corollary~\ref{char_binom_compl_inter_cor} we show that a necessary and sufficient condition for $I_L$ to be a binomial complete intersection ideal, is  either for $L$ to be equal to $L_{pure}$ or for a basis of $L^\sigma$ to have vectors giving the rows of a mixed dominating matrix.  Thus Corollary~\ref{char_binom_compl_inter_cor}  
completely characterizes binomial complete intersection lattice ideals for arbitrary lattices. 

In Section 6 we work in detail an example: we use the techniques developed in this paper, in order to compute all cardinality-minimal generating sets of a lattice ideal.

\section{Fibers, the Pure Sublattice and Bases of a Lattice}

Let $L$ be a lattice in $\mathbb{Z}^n$, $R=\Bbbk[x_1,\ldots, x_n]$ be the polynomial ring over a field $\Bbbk$ and   $I_L\subseteq R$ be the lattice ideal:
\[
I_L=\langle x^{\ub}-x^{\vb}:\ \ub,\vb\in\N^n, \ {\ub}-{\vb}\in L\rangle. \]
\begin{def1}\label{mu_I_L}
{\em We let $\mu(I_L)$ be the minimal cardinality of a minimal generating set of $I_L$. In other words, $\mu(I_L)$ is the cardinal of a cardinality-minimal generating set of $I_L$.}
\end{def1}

We denote by $\mathbb{T}^n$ the set of monomials of $R$ including $1=x^{\bf 0}$ and by $\N$ the set $\Z_{\geq 0}$. If $J$ is a monomial ideal of $R$ we denote  by $G(J)$ the unique minimal set of  monomial generators of $J$. For $r\in \N$ we let $[r]=\{1,\ldots, r\}$. Let $\ab=(a_1,\ldots, a_n)$, $\bb=(b_1,\ldots, b_n)\in \Z^n$. We write $\ab\geq{\bf 0}$ if $\ab\in \N^n$. If either $\ab\geq{\bf 0}$ or $-\ab\geq{\bf 0}$ we say that $\ab$ is {\it pure}. With this terminology, a lattice $L$ is non-positive if and only if  it contains a nonzero pure vector, otherwise  it is positive. We write $\ab\geq\bb $ if $a_i\geq b_i$ for $i=1,\ldots, n$, and we say that $\ab$, $\bb$ are incomparable if $\ab-\bb$ is not pure. In general we let $\supp(\ab)=\{i: a_i\neq 0\}\subset [n]$. For any subset $X$ of $\Z^n$ we  let
\[
\supp(X):=\bigcup_{w\in X} \supp(w) \ .
\]

\begin{def1}\label{fiberdef}
\rm{ We say that $F\subset \mathbb{T}^n$ is an $I_L$-{\it  fiber} if there exists  $x^{\ub}\in \mathbb{T}^n$ such that $F= \{ x^{\vb}\in \mathbb{T}^n:\ {\vb}-{\ub}\in L \}$. If $x^{\ub}\in F$, and  $F$ is an $I_L$-fiber we write $F_{\ub}$ or $F_{x^{\ub}}$ for $F$. If  $B\in I_L$ and $ B=x^{\ub}-x^{\vb}$ we  write $F_B$ for $F_{\ub}$. When $F$ is an $I_L$-fiber we let $M_F=\langle x^{\ub}: \ x^{\ub}\in F \rangle$ be the monomial ideal generated by the elements of $F$.}
\end{def1}

\noindent From the properties of the lattice and the definition of lattice ideals we get the following:

\begin{prop1} If $x^{\vb}\in F_{\ub}$ then $F_{\ub}=F_{\vb}$. Moreover $F_{\ub}=\{x^{\vb}: x^{\vb}-x^{\ub}\in I_L\}$. If $x^{\ub}-x^{\vb}\in I_L$ then ${\ub}-{\vb}\in L$.
\end{prop1}

\noindent We  remark that $F_{\ub}$ is a singleton if and only if there is no binomial $0\neq B\in I_L$ such that $F_B=F_{\ub}$. We note that $F\subset M_F$ and  $G(M_F)\subset F$. The following proposition follows also from \cite[Theorem 8.6]{MS}.

\begin{prop1}\label{infinitefiber} 
Let $L\subset \Z^n$ be a lattice. The following are equivalent:
\begin{enumerate}
 \item  The lattice $L$ is non-positive.
\item All $I_L$-fibers are infinite.
\item There exists an $I_L$-fiber which is infinite.
\end{enumerate}
\end{prop1}

\begin{proof} 
$(1)\Rightarrow (2):$ Let $F$ be an $I_L$-fiber and suppose that ${\bf 0}\neq {\ub}\in L\cap \N^n$. It is easy to see that if ${\vb}\in \N^n$ and $F$ is the $I_L$-fiber such that $x^{\vb}\in F$, then $x^{\vb+l\ub}\in F$ for all $l\in \N$. Thus $F$ is infinite. Implication $(2)\Rightarrow (3)$ is obvious. For $(3)\Rightarrow (1)$ suppose that an $I_L$-fiber $F$ is infinite. Let $x^{\vb}\in F$ be such that $x^{\vb}\notin G(M_F)$. Note that since $F$ is infinite such a $\vb$ exists. Since $x^{\vb}\in M_F$, there exists a monomial $x^{\ub}\in G(M_F)$ such that $x^{\ub}|x^{\vb}$ and thus $x^{\vb}=x^{\wb}x^{\ub}$ for ${\bf 0}\neq {\wb}\in \N^n$. Since $x^{\vb},x^{\ub}\in F$ it follows that $\wb=\vb-\ub\in  L$. Therefore  ${\bf 0}\neq\wb\in L\cap \N^n$ and consequently $L$ is non-positive.
\end{proof}

\begin{cor1}\label{gen_fin_fiber_cor} 
Let $L\subset \Z^n$ be a lattice.  The lattice $L$ is positive if and only if $G(M_F)=F$, where $F$ is any $I_L$-fiber.
\end{cor1}

\begin{proof}
Suppose that $L\cap \N^n=\{\bf 0\}$. Since $M_F=\langle F\rangle$ to prove that $G(M_F)=F$, it is enough to show that if $x^{\ab}\neq x^{\bb}\in F$  then $\ab$ and $\bb$ are incomparable. Suppose otherwise. Then $\ab-\bb\in L$ is pure, a contradiction. For the other direction suppose that $G(M_F)=F$. Thus $F$ is finite and the conclusion follows from Proposition \ref{infinitefiber}.
\end{proof}

In the following definition we introduce  the induced lattices that appear in our work and  the related notation. 

\begin{def1}\label{Lpure}
{\rm We let $L^+=L\cap \N^n$, $\sigma=\supp(L^+)$ and $L_{pure}$ be the subgroup of $L$ generated by $L^+$. If $\ub=(u_i)\in \Z^n$ and 
$\sigma=[n]$, then we set $\ub^{[n]}=0$. Otherwise we let $\ub^{\sigma}$  be the vector $(u_i)_{i\notin \sigma}$. We define $L^{\sigma}$ to be the lattice  generated by the vectors $\ub^\sigma$ where $\ub\in L$. Finally, if $\sigma\neq \emptyset$, we let   $\ub_{\sigma}=(u_i)_{i\in\sigma}$ and $(L_{pure})_\sigma$ to be the lattice generated by  the vectors $\ub_{\sigma}$,   $\ub\in L_{pure}$. If $\ub=(u_i)\in \Z^n$ and 
$\sigma=\emptyset$, then we set $\ub_{\emptyset}=0$.}
\end{def1}

In the course of the proof of Proposition \ref{infinitefiber} we proved the following:

\begin{prop1}\label{descfiber} 
Let $L\subset \Z^n$ be a lattice and let $F$ be an $I_L$-fiber. If $G(M_F)=\{x^{\ab_1},\ldots,x^{\ab_s}\}$ then 
\[
F=\bigcup^s_{i=1}\{x^{\ab_i}x^{\wb}:\ \wb\in L^{+}\}\ .
\]
\end{prop1}

The next proposition considers the support of the elements of $L$ that belong to $L_{pure}$. 

\begin{prop1}\label{supp_pure_prop} 
There exists an element $\wb$ in $L^+$ such that $\supp(\wb)= \sigma$. For $\ub\in L$ we have that $\supp(\ub)\subset \sigma$ if and only if $\ub\in L_{pure}$.
\end{prop1}
\begin{proof}
The existence of $\wb$ follows from the observation that if $\wb_1,\wb_2\in L^+$ then $\wb_1+\wb_2\in L^+$   and $\supp(\wb_1)\cup \supp(\wb_2)=\supp(\wb_1+\wb_2)$. 
 
Suppose now that $\ub\in L$ and $\supp(\ub)\subset \sigma$. Let $\wb\in L^+$ be such that $\supp(\wb)=\sigma$. It is clear that for $l\in \N$, $l\gg 0$, $\ub+l\wb= \wb'\in \N^n$. Since $\ub$, $l\wb\in L$ it follows that $\wb'\in L$ and thus $\wb'\in L^+$. Therefore    $\ub=\wb'-l\wb\in L_{pure}$. 
\end{proof}

The following is an immediate consequence of Proposition~\ref{supp_pure_prop}.

\begin{cor1}\label{pure_supp_outside} 
Let $L\subset \Z^n$ be a lattice and $\ub\in L$. Then $\ub\in L_{pure}$ if and only if $\ub^{\sigma}={\bf 0}$.
\end{cor1}

We note that since $L_{pure}$ is generated by the elements of $L^+$ it is clear that 
\[
 \supp(L_{pure})=\sigma \ .
\]

\begin{def1}{\rm A nonzero vector $\ub\in L$ is called $L$--primitive if whenever $\lambda \ub\in L$ where $\lambda \in  \Q$ then $\lambda\in \Z$. In other words, $\ub$ is $L$--primitive if and only if $\Q \ub\cap L=\Z \ub$. }
\end{def1}

\noindent Equivalently $\ub$ is $L$--primitive if it is the ``smallest" element of $L$ in the direction determined by $\ub$.

\begin{prop1}\label{prim_dir}
Let ${\bf 0}\neq \vb\in L$. There is an $L$--primitive vector $\ub\in L$ such that $\vb=\lambda \ub$ for $\lambda\in \Z$.
\end{prop1}

\begin{proof} Since $\vb\in L\setminus\{\bf 0\}$ it follows that $\Q\vb\cap L$ is a nonzero subgroup of $L$. Hence $\Q\vb\cap L$ is a free abelian group of rank $1$, whose generator is the desired $\ub$.
\end{proof}

\noindent  Consider now any basis  of $L$ as a $\Z$-module. The next theorem states that the elements of such a basis are necessarily $L$--primitive.

\begin{thm1}\label{primitive_bases} 
Let  $L\subset \Z^n$ be a lattice and let  $\B$ be a basis of $L$ as a $\Z$-module. The elements of $\B$ are $L$--primitive.
\end{thm1}

\begin{proof} Since $L$ is a sublattice of $\Z^n$, there exists an $r\in \N$ such that $L\cong \Z^r$, $r=\rank (L)$.  Let  $\B=\{\ub_1,\ldots, \ub_r\}$. Since $\B$ is a basis then it follows immediately that $\Q\ub_i\cap L=\Z\ub_i$ for all $i$. 
\end{proof}

The next theorem shows that given an $L$--primitive vector, one can find a basis of 
$L$ that contains it.
\begin{thm1}\label{basis_primitive_lattice} 
Let $L$ be a lattice and $\ub$ an $L$--primitive vector. There exists a basis $\B$ of $L$ such that $\ub\in\B$.
\end{thm1}

\begin{proof} Since $L$ is a lattice then $L/{\Z\ub}$ is a finitely generated abelian group. If $L/{\Z\ub}$ is not torsion-free then there exists a $\vb\in L\setminus \Z\ub$ such that $\lambda\vb\in \Z\ub$ for some positive integer $\lambda>1$. This implies that $\Q\ub\cap L\supsetneq\Z\ub$, a contradiction to the fact that $\ub$ is $L$-primitive. Thus $L/{\Z\ub}$ is torsion-free and consequently a free abelian group of rank equal to $\rank(L)-1$. Now lifting the elements of any basis of $L/{\Z\ub}$ to 
elements of $L$ and adding $\ub$ we obtain a basis of $L$.
\end{proof}

 As the next corollary shows, we can find bases of $L_{pure}$ whose elements are in $\N^n$  having full support $\sigma$.

\begin{cor1}\label{cor_basis_pure}
Let $L$ be a non-positive lattice. There exists a basis of $L_{pure}$ whose elements are in $L^+$ and have support equal to $\sigma$.
\end{cor1}

\begin{proof} By Proposition~\ref{supp_pure_prop} and Proposition~\ref{prim_dir} there is an $L$-primitive  vector $\ub_1\in L^+$  such that $\supp(\ub_1)=\sigma$. By Theorem~\ref{basis_primitive_lattice} there exists a basis $\{\ub_1,\ldots,\ub_r\}$ of $L_{pure}$. It is clear that for $l\gg 0$, $\ub'_i=\ub_i+l \ub_1\in L^+ $ for $i=2,\ldots, r$.  The set $\{ \ub_1, \ub_2',\ldots, \ub_r'\}$ has the desired properties. 
\end{proof}

Bases of the lattice $L$ are clearly important for the study of $I_L$, even though they do not directly give the elements of the minimal generating sets of $I_L$. Indeed, it is well known that the following holds:
\[ 
I_L= (x^{\ub_1^+}-x^{\ub_1^-},\ldots,x^{\ub_r^+}-x^{\ub_r^-}):\ (x_1\cdots x_n)^\infty,
\] 
where $\ub_1,\ldots,\ub_r$ form a basis of  $L$ and 
${\bf u_i}^+,{\bf u_i}^-\in\N^n$ are (for all $i$) the unique vectors of disjoint support,
such that ${\bf u_i}={\bf u_i}^+-{\bf u_i}^-$, see  \cite[Lemma 12.2]{St}. 

In the next two sections we describe all minimal generating sets of lattice ideals using properties of the $I_L$-fibers and of the bases of $L_{pure}$.

\section{Fibers and minimal generating sets of Lattice Ideals}

 Let  $L\subset \Z^n$  be a lattice and $R=\Bbbk[x_1,\ldots, x_n]$, where $\Bbbk$ is a field. If $G\subset \mathbb{T}^n$ and $\vb\in \N^n$ we let $x^{\vb} G:=\{ x^{\vb} x^{\ub}:\  x^{\ub}\in G\}$.

\begin{lem1}\label{prop_fiber_mult}
Let $G$, $F$ be   $I_L$-fibers. If there exists  $\wb_1\in \N^n$, $x^{\ub}\in G$  such that  $x^{\wb_1}x^{\ub}\in F$  then $x^{\wb_1}G\subset F$. Moreover if $x^{\wb_2}F\subset G$ for $\wb_2\in \N^n$ then $\wb_1+\wb_2\in L^+$ and $\supp(\wb_i)\subset \sigma$, $i=1,2$.
\end{lem1}
\begin{proof} Suppose that $x^{\wb_1} x^{\ub}=x^{\vb}\in F$ and let $x^{\ub'}\in G$. Since $\ub'-\ub\in L$ it follows that $(\wb_1+\ub')-\vb \in L$ and thus $x^{\wb_1+\ub'}\in F$.

 If in addition $x^{\wb_2}F\subset G$ it follows that $x^{\wb_1+\wb_2}F\subset F$. Since $x^{\wb_1+\wb_2}x^{\vb}$, $x^{\vb}\in F$ it follows that $\wb_1+\wb_2\in L\cap\N^n $.
 \end{proof}

\begin{def1}\label{equiv} 
{\em Let  $F$, $G$ be $I_L$-fibers. We say that $F\equiv_{_{I_L}}G$ if there exist $\ub,\vb\in \N^n$ such that $x^{\ub}F\subset G$ and $x^{\vb} G\subset F$.}
\end{def1}
\noindent It is immediate that  $F\equiv_{_{I_L}}G$ is  an equivalence relation among the $I_L$-fibers. We denote the equivalence class of $F$ by $\overline F$. Thus
\[ 
\overline{F}=\{ G:\  G \textrm{ is an } I_L\ \textrm{-fiber},\ G \equiv_{_{I_L}}F\}\ .
\]
We note that $F\equiv_{_{I_L}}G$ implies that the cardinality of $F$ is equal to the cardinality of $G$.

\begin{lem1}\label{finite_equivalence_fiber} 
If $L$ is positive and $F$ is an $I_L$-fiber then $\overline F=\{ F\}$.
\end{lem1}

\begin{proof} By Proposition \ref{infinitefiber}, $|F|<\infty$. Let $G$ be an $I_L$-fiber, $G\equiv_{_{I_L}}F$. There are $\ub,\vb\in \mathbb{N}^n$ such that $x^{\ub} F\subset G$ and $x^{\vb}G\subset F$. Since $|F|=|x^{\ub}F|=|G|=|x^{\vb}G|$ it follows that   $x^{\vb} x^{\ub} F= F$ and $x^{\vb}=x^{\ub}=1$.
\end{proof}

Next we want to investigate the number of equivalent fibers inside each equivalence class when $L_{pure}\neq\{\bf 0\}$ and  $\sigma\neq \emptyset$.  First we note the following:

\begin{lem1}\label{equiv_rema}
{\em Let $\sigma\neq \emptyset$, $F$ an $I_L$-fiber and $\ub\in \N^n$ such that $\supp(\ub)\subset \sigma$. If $G$ is an $I_L$-fiber with the property $x^{\ub} F\subset G$ then $G\in \overline F$.}
\end{lem1}
\begin{proof} Let $\wb\in L^+$ be such that $\supp(\wb)=\sigma$. There exists $l\gg 0$ such that $l\wb-\ub\in \N^n $. Since $l\wb\in L$ it follows that  $x^{l\wb}F\subset F$. Let $x^{\vb}\in G$ such that $x^{\vb}=x^{\ub} x^{\bf p}$ for $x^{\bf p}\in F$. It follows that $x^{l\wb-\ub}x^{\vb}=x^{l\wb} x^{\bf p}\in F$ and thus $x^{l\wb-\ub} G\subset F$ by Lemma~\ref{prop_fiber_mult}.
\end{proof}

We recall that if $\ub\in \Z^n$,  then   $\ub_{\sigma}=(u_i)_{i\in\sigma}$. 
Thus, if $s=|\sigma|$, then we can assume that $\ub_{\sigma}\subset \Z^s$ and that
$(L_{pure})_\sigma$ is a sublattice of $\Z^s$.

\begin{prop1} \label{same_nr_equiv} 
Let $L$ be a non-positive lattice and $F$ an $I_L$-fiber. The cardinality of  $\overline{F}$ is equal to $|\Z^s/(L_{pure})_\sigma|$, where $s=|\sigma|$.
\end{prop1}
\begin{proof} For every $G\in \overline F$ choose $\ub_G\in \N^n $ such that $x^{\ub_G} F\subset G$. Let 
\[ 
\phi:\overline F\to \Z^s/(L_{pure})_\sigma,\ \ \phi(G)=(\ub_G)_{\sigma} + (L_{pure})_\sigma\ .
\]
The definition of $\phi$ is independent of the choice of $\ub_G$. Indeed suppose that $\ub,\vb\in \N^n$ are such that $x^{\ub} F\subset G$ and $ x^{\vb} F\subset G$. This implies that $\ub-\vb\in L$. By Lemma~\ref{prop_fiber_mult} it follows that $\ub^{\sigma}=\vb^{\sigma}={\bf 0}$ and by Proposition \ref{supp_pure_prop}   we obtain $\ub-\vb\in L_{pure}$ and $\ub_{\sigma}-\vb_{\sigma}\in (L_{pure})_\sigma$.

We will show that   $\varphi$ is a bijection: the only part needing proof is the surjectivity of $\varphi$.  Let $\ub'+ (L_{pure})_\sigma$ be an element of  $\Z^s/ (L_{pure})_\sigma$.  First we remark that we can assume without loss of generality that $\ub'\in \N^s$.  Indeed, let   $\wb\in L^+$ be such that $|\supp(\wb_\sigma)|=s$. It is clear that for $l\gg 0$, $l \wb_\sigma +\ub'\in \N^s$ and thus  $\ub'+(L_{pure})_\sigma=(l\wb+\ub')+(L_{pure})_\sigma$. Let $\ub\in\N^n$ be such that $\ub^{\sigma}={\bf 0}$, $\ub_{\sigma}=\ub'$ and   let $G$ be the  $I_L$-fiber such that $x^{\ub}F\subset G$. By Remark \ref{equiv_rema} it follows that $G\in \overline F$, and thus $\phi(G)= \ub'+ (L_{pure})_\sigma$.
\end{proof}

\begin{ex2}\label{example_multi} 
{\rm (a) We consider the lattice ideal $I_L=\langle 1-xy\rangle \subset \Bbbk[x,y]$ where  $L=\langle (1,1) \rangle\subset \Z^2$.  There are infinitely many  $I_L$-fibers: for any $c\in \mathbb{Z}$ the set $F_c=\{x^iy^j:\ i-j=c\}$ is an $I_L$-fiber. All $I_L$-fibers are infinite and belong to the same equivalence class: the cardinality of this equivalence class is  $|\Z|$. Indeed    $L_{pure}=L$,  and $\Z^2/\ L\cong \Z$.

(b) If we consider the lattice ideal $I_L=\langle 1-xy,1-x^5\rangle\subset\Bbbk[x,y]$, where $L=\langle (1,1),(5,0)\rangle\subset\Z^2$ then there are exactly five infinite $I_L$-fibers:
\[
F_{k}=\{x^iy^j\ :\ i-j\equiv k\mod 5\}, \quad 0\leq k\leq 4\ , 
\]
which are all equivalent. Hence we have only one equivalence class $\overline {F_0}=\{F_0,\ldots,F_4\}$ which has five equivalent fibers. Indeed  $L_{pure}=L$ and $\Z^2/L\cong\Z_5$.}
\end{ex2}

\noindent We define the relation ``$\leq_{_{I_L}}$" among the equivalence classes of $I_L$-fibers.

\begin{def1}\label{order} 
{\em Let $F$, $G$ be $I_L$-fibers. We say that  $\overline{F}\leq_{_{I_L}} \overline{G}$ if there exists $\ub \in \N^n$ such that $x^{\ub} F \subset G$.} 
\end{def1}

\noindent It is immediate that   ``$\leq_{_{I_L}}$" is well defined and is a partial order among the equivalence classes of $I_L$-fibers. For simplicity of notation we occasionally  write $F\leq_{_{I_L}} G$ if $\overline F \leq_{_{I_L}} \overline G$ and   $F<_{_{I_L}} G$ if $\overline F\leq_{_{I_L}} \overline G$ and $\overline F\neq \overline G$. We  note that  $F_{\{1\}}\leq_{_{I_L}} F$ for any $I_L$-fiber $F$. We also remark that if $L$ is positive then  $<_{_{I_L}}$ gives the ordering on the fibers of $I_L$ induced by the $\Z^n/L$-degrees, see \cite[Section 2]{CKT}.

\begin{thm1}\label{chain_least} 
Any strictly descending chain of  equivalence classes of $I_L$-fibers is finite.
\end{thm1}

\begin{proof} Assume that
\[ 
 \overline F_1\ >_{_{I_L}}\ \cdots\  >_{_{I_L}}\   \overline F_k\ >_{_{I_L}} \ \overline F_{k+1}\  >_{_{I_L}}\  \cdots 
\]
is a chain of equivalence classes of fibers with no least element.  Choose a representative $F_i$, $i\in \N$ for each class. Next  consider the corresponding ascending chain of monomial ideals: 
\[
 M_{F_1} \subset \cdots \subset M_{F_1}+\cdots+M_{F_k}\subset M_{F_1}+\cdots+M_{F_{k+1}}\subset \cdots
\]
The chain stabilizes at some step, say $s$, so that 
\[ 
M_{F_1}+\cdots+M_{F_s}= M_{F_1}+\cdots+M_{F_{s+1} }.
\]
Let $x^{\ab}\in G(M_{F_{s+1}})$. By the above equality it follows that $x^{\ab}\in M_{F_i}$ for some $1\le i<s+1$ and $x^{\ab}=x^{\ub} x^{\bb}$ where $x^{\bb}\in G(M_{F_i})$. Since $x^{\ub}x^{\bb}\in F_{s+1}$ it follows that $x^{\ub} F_i \subset F_{s+1}$. This leads to a contradiction since $\overline{F}_{s+1}<_{_{I_L}} \overline{F}_i$.
\end{proof}
\noindent


\begin{def1}\label{def_markov_fiber}
\rm{ Let $F$ be an $I_L$-fiber. We say that $F$ is a {\it Markov fiber} if there exists a cardinality-minimal generating set $S$ of $I_L$ such that $\overline F=\overline F_B$ for some $B$ in $S$.}
\end{def1}

Next, we generalize the constructions of  \cite[Section 2]{CKT} to arbitrary lattice ideals. Let $F$ be an $I_L$-fiber. We let
\[
 I_{{L,{<\overline F}}}=\langle B\in I_L:\ B\ \rm{ binomial},\ \overline F_B <_{_{I_L}} \overline F  \rangle 
\]
and
\[
I_{L,{\leq \overline F}}=\langle B\in I_L:\ B\ \rm{ binomial},\ \overline F_B\leq_{_{I_L}} \overline F \rangle.
\]
We note that  $I_{L,{<\overline F}}=0$ if there is no $I_L$-fiber $G$ such that $\overline G<_{_{I_L}} \overline F$. It is  clear that the definition of these ideals does not depend on the chosen fiber representative. Finally  if $S$ is any subset of binomials of $I_L$ we let
\[ 
S_{\overline F}=\{ B\in S: \overline F_B\leq_{_{I_L}} \overline F\}. 
\]

\begin{rem1}\label{fiber_of_one}{\rm  We will pay extra attention to the fiber that contains $1$, $F_{\{1\}}$. Let $S$ be a set of binomials of $I_L$. According to the definitions \[ S_{\overline F_{\{1\}}}=\{B\in S:\ F_B\in  \overline F_{\{1\}}\} \textrm{ and } I_{L_{pure}}=I_{L,{\le \overline F_{\{1\}}}}.\] }
\end{rem1}

We isolate the following  proposition whose proof is within the
proof of \cite[Lemma A.1]{HM}.

\begin{prop1}\label{sum_generators} Let $S$ be a minimal generating set of $I_L$, $F$ an $I_L$-fiber, $x^{\wb_1}$, $x^{\wb_2}\in F$. There exists a subset $T\subset S_{\overline F}$ such that 
\[ 
x^{\wb_1}-x^{\wb_2}= \sum_{i,B} \pm x^{\ab_{i,B}} B, 
\]
where $B\in T$, $\ab_{i,B}\in \N^n$ and  $\ab_{i,B}\neq \ab_{j,B}$ for $i\neq j$.
 \end{prop1}

\noindent We note that in the summation, the binomials $B\in T$ may 
appear more than once.  The emphasis of the above statement is  that when summing up and factoring out the binomials $B\in T$, we get an expression 
\[
 x^{\wb_1}-x^{\wb_2}= \sum f_i B_i  
\]
where $B_i\neq B_j\in T$, for $i\neq j$, and all nonzero coefficients of the monomial  
terms of the $f_i$ are  $\pm 1$. 
Next we describe the ideals $I_{L,{< \overline F}} $ and $I_{L,{\le \overline F}}$ in terms of the generators of $I_L$.

\begin{prop1} \label{leq_prop}
Let  $S$ be a generating system of binomials for $I_L$. The following hold:
\[
 I_{_{L,{<\overline F}}}=\langle B:\  B\in S,\  \overline F_{B}<_{_{I_L}} \overline F \rangle
\] 
and
\[ 
I_{L,{\leq \overline F}}=\langle B:\  B\in S, \ {\overline F}_B\le_{_{I_L}}{\overline F} \rangle.
\]
\end{prop1}

\begin{proof} We will  show the statement for $I_{L,{<\overline F}}$, the other one having a similar proof.  Let $J=\langle B:\  B\in S,\  \overline F_{B}<_{_{I_L}} \overline F \rangle$. We will show that $J=I_{L,{<\overline F}}$. It is clear that  $J\subset I_{L,{< \overline F}}$. To show the other containment it is enough to show that if  $B=x^{\ub}-x^{\vb}\in I_{L,{<\overline F}}$ then $B\in J$. Let $B=x^{\ub}-x^{\vb}\in I_{L,{<\overline F}}$. Since $B\in I_L$, by Proposition \ref{sum_generators} it follows that $B=\sum_{i=1}^t \pm x^{\ab_{i,B_i}}B_{i}$ where $B_{i}\in S$ are not necessarily distinct while $\ab_{i,B_i}\neq \ab_{j,B_j}$ for $B_i=B_j$ and $i\neq j$. We will do induction on $t$.  Without loss of generality we can assume that $B_{1}=x^{\ub_1}-x^{\vb_1}$ and $x^{\ab_{1}}x^{\ub_1}=x^{\ub}$, the other cases being done similarly. First we show the inductive step. Suppose that $t=1$. Since $x^{\ab_{1}}x^{\ub_1}=x^{\ub}$ it follows that $x^{\ab_{1}}  F_{B_1}\subset F_B$. Thus $\overline F_{B_1}\le_{_{I_L}} \overline F_B$. Since $\overline  F_B<_{_{I_L}}\overline  F$ we see that $\overline F_{B_1}<_{_{I_L}} \overline F$. Assume now that $t>1$ and consider $B'=B-x^{\ab_{1}}B_1=x^{\ab_{1}}x^{\vb_1}-x^{\vb}$. Since $F_{B'}=F_B$, it follows that $B'\in I_{L,{<\overline F}}$ and we are done by induction.
\end{proof}

We can now determine when the $I_L$-fiber $F$ is a Markov fiber, in terms of the subideals of $I_L$, defined earlier.

\begin{thm1} \label{weak_betti} 
Let $F$ be an $I_L$-fiber. $F$ is a Markov fiber if and only if $I_{L,{< \overline F}}\not= I_{L,{\leq \overline F}}$.
\end{thm1}

\begin{proof} Let $S$ be a cardinality-minimal generating set of $I_L$. If $I_{{L,{<\overline F}}}\neq I_{\leq \overline F}$ then by Proposition~\ref{leq_prop} there exists a $B\in S$ such that $F_B \equiv_{_{I_L}} F$. For the converse assume that there exists $B\in S$ such that $\overline F_B=\overline F$. It follows immediately that  $I_{{L,{<\overline F}}}\neq I_{L,{\leq \overline F}}$.
\end{proof}

\begin{cor1}\label{markov_fib_inv_cor} 
The set of equivalence classes of Markov fibers of a lattice ideal $I_L$ is an invariant of $I_L$.
\end{cor1}

In the following section, we will refine this result taking under account the number of times an equivalence class of a Markov fiber appears in a cardinality-minimal generating set of $I_L$, see Theorem~\ref{cor_multiset}.

\section{Generating sets of lattice ideals}

\bigskip
\subsection{Pure Lattices} First we consider the case of $L\subset \Z^n$ being a pure lattice, i.e.~when $L$ is generated by its pure elements and show how to obtain all cardinality-minimal generating sets of $I_L$. Let $S=\{B_1,\ldots, B_r\}$ be a set of binomials of $I_L$. We say that $S'$ is a {\it rearrangement} of $S$ if there is a bijective function $f: S\rightarrow S'$ such that $f(B_i)=\pm B_j$. Compositions of rearrangements are rearrangements. It is clear that if $S$ is a generating set of $I_L$ then all rearrangements of $S$ are generating sets of $I_L$. The theorem below generalizes \cite[Lemma 2.1]{SWZ}.

\begin{thm1}\label{Markov_pure} Let  $L=L_{pure}$ be a  pure lattice of rank $r$, $\sigma=\supp(L)$ and $S$ a set of $r$ binomials of $I_L$. The set $S$
generates $I_L$ if and only if there is a rearrangement $\{ x^{\ub_1}-x^{\vb_1}, x^{\ub_2}-x^{\vb_2},\ldots, x^{\ub_r}-x^{\vb_r}\}$ of $S$ such that the following three conditions are satisfied:
\begin{enumerate}
\item{}$\{\ub_1-\vb_1, \dots , \ub_r-\vb_r\}$ is a basis of $L$,
\item{} for $i\in [r]$, $\supp(\ub_i)\cup \supp(\vb_i)\subset \sigma$,
\item{} $(a)$ $x^{\vb_1}=1$ and $(b)$ $\supp (\vb_i)\subset \bigcup_{j=1}^{i-1}\supp (\ub_j)$ for $2\leq i\leq r$.
\end{enumerate}
\end{thm1}
\begin{proof} Suppose first that $S=\{B_1, \dots,B_r\}$ generates $I_L$. We let $B_j= x^{\bb_j}-x^{\cb_j}$ for $j\in [r]$. Let $\ub\in L$. According to Proposition~\ref{sum_generators} there is an index set $A$ such that $$x^{\ub^+}-x^{\ub^-}=\sum_{l\in A} \pm x^{\ab_l}(x^{\bb_j}-x^{\cb_j}),$$ where $i_l\in [r]$, $\ab_l\in \N^n$. Expanding the RHS, equating the exponents of the equal monomial terms, subtracting the expressions for $\ub^+$ and $\ub^-$ and substituting  $\ab_l$, $l\in A$,  one gets that $\ub\in \Z (\bb_1-\cb_1)+\cdots+\Z (\bb_r-\cb_r)$. This shows that $L=\Z (\bb_1-\cb_1)+\cdots+\Z (\bb_r-\cb_r)$.  Since $\rank (L)=r$  it follows that $\{\bb_1-\cb_1, \dots , \bb_r-\cb_r\}$ is a basis of $L$. Next we remark that  $\supp (\bb_i) \subset \sigma$ if and only if $\supp(\cb_i)\subset \sigma$.  Indeed this is immediate since $\supp(\bb_i-\cb_i)\subset \sigma$. Now suppose that for some   $i\in [r]$, $\supp (\bb_i)\not \subset \sigma$. We claim that $B_i=x^{\bb_i}-x^{\cb_i}$ is redundant in $S$  as a generator of $I_L$. For this we will show that if $\ub\in L$ then $x^{\ub^+}-x^{\ub^-}$ can be written as a linear combination of the  elements of  $S\setminus\{B_i\}$. Indeed consider again the relation 
\[
x^{\ub^+}-x^{\ub^-}=\sum_{l} \pm x^{\ab_l}(B_{j_l})
\] 
of Proposition~\ref{sum_generators}. Substitute the value $0$ to any variable $x_j$ where $j\notin \sigma$: the terms in the above relation involving $B_i$ disappear. Thus  $S\setminus \{B_i\}$ is a generating set of $I_L$. This is of course a contradiction by  the generalized Krull's Principal Ideal Theorem, since the height of $I_L$ is $r$, see \cite[Theorem 2.1]{ES}. To show that there is a rearrangement of $S$ that satisfies the conditions of the theorem we notice that $S$ contains a binomial $B_j$ with 1 as one of its monomial terms. Indeed let $\wb\in L^+$ with $\supp(\wb)=\sigma$, see Proposition~\ref{supp_pure_prop}. Since  $x^{\wb}-1\in I_L$, $x^{\wb}-1=\sum \pm x^{\ab_l}(B_{i_l}).$  It is clear that there exists a value of $l$ such that  a monomial term of $B_{i_l}$ is equal to 1 (and $\ab_{l}=0$):  otherwise $x^{\wb}-1\in \langle x_1,\dots ,x_n\rangle$,  a contradiction. It is immediate  that we  can rearrange $S$ by a bijective function
$f_1$ so that $f_1(B_{i_l})=x^{\ub_1}-1 $. Next we claim that there is  $B=x^{\bb}-x^{\cb}\in S$ such that $B\neq B_{i_l}$ and $\supp(\bb)$ or $\supp(\cb)\subset \supp(\ub_1)$. Indeed, suppose not. Then clearly $\supp(\ub_1)\neq \sigma$. Consider again the expression 
\[
x^{\wb}-1=\sum_{t\in A, i_t=i_l} x^{\ab_t} B_{i_l}+\sum_{t\in A, i_t\neq i_l} \pm x^{\ab_t}(B_{i_t})\ .
\]
Substitute the value 1 for all variables whose index is in $\supp(\ub_1)$ and the value $0$ for all other variables. We obtain a contradiction: $-1=0$. To avoid the contradiction there must be $B\in S$ so that $\pm B=x^{\ub_2}-x^{\vb_2}$ and $\supp(\vb_2)\subset \supp(\ub_1)$. We rearrange $f_1(S)$ by $f_2$ which keeps all elements of $f_1(S)$ fixed but $B$: $f_2(B)=x^{\ub_2}-x^{\vb_2}$.  More generally once $f_s$ has been defined for $s<r$ so that the third condition is satisfied for all $i \le s$, the same argument produces $f_{s+1}$ with the desired property.

We now prove the converse. Consider a set $S$ of binomials whose rearrangement $\{ x^{\ub_1}-1, x^{\ub_2}-x^{\vb_2},\ldots, x^{\ub_r}-x^{\vb_r}\}$ satisfies  the three conditions of the theorem. Let $J=\langle x^{\ub_1}-1, \dots,x^{\ub_r}-x^{\vb_r} \rangle$. We will show that $J=I_L$. It is clear that $J\subset I_L$. Since $\ub_1,\dots, \ub_r-\vb_r$ is a basis of $L$ and $\bigcup_{i=1}^r(\supp (\ub_i)\cup \supp(\vb_i)) \subset \sigma$ it is clear that $\bigcup_{i=1}^r(\supp(\ub_i)\cup \supp(\vb_i)) = \sigma$. By the third condition it follows that 
\[
\bigcup_{i=1}^r \ \supp (\ub_i) =\sigma\ .
\]
Next we will show  that for every $k\in [r]$ there exists $\wb_k\in L^+$ such that $x^{\wb_k}-1\in J$ and $\supp(\wb_k)=\bigcup_{j=1}^{k}\supp (\ub_j)$. For $k=1$ we set $\wb_1=\ub_1$.  Since $\supp(\vb_2)\subset \supp(\ub_1)$ there exists $\lambda_1\in \N$, $\lambda_1\gg 0$ such that $\lambda_1 {\wb}_1>{\vb}_2$. We set ${\wb}_2=( \lambda_1 {\wb}_1-\vb_2)+\ub_2$; $\supp(\wb_2)=\supp(\ub_1)\cup \supp(\ub_2)$. Moreover 
\[
x^{\wb_2}-1=  x^{\lambda_1 \wb_1-\vb_{2}}(x^{\ub_{2}}-x^{\vb_{2}})+x^{\lambda_1 \wb_1}-1\in J\ ,
\]
as wanted. It is clear that this construction generalizes for all $k\in [r]$. In particular $\supp(\wb_r)=\sigma$.

We will now show that if $\ub-\vb\in L$ then $x^{\ub}-x^{\vb}\in J$, which ends the proof. Since $J:\ (x_1\cdots x_n)^\infty=I_L$ there exists $\wb\in \N^n$ such that $x^{\wb}(x^{\ub}-x^{\vb})\in J$. By the second condition it is clear that $\wb$ can be chosen so that $\supp(\wb)\subset \sigma$. Since $\supp(\wb_r)=\sigma$ there exists $\lambda\in \N$, $\lambda\gg 0$ such that $\lambda \wb_r>\wb$. Therefore $x^{\lambda \wb_r} (x^{\ub}-x^{\vb})\in J$. It follows that 
\[
 x^{\ub}-x^{\vb}= (x^{\lambda \wb_r}-1) (x^{\vb}-x^{\ub})-x^{\lambda \wb_r} (x^{\ub}-x^{\vb})\in J 
\]
and consequently  $I_L=J$.
\end{proof}

We remark that the binomials of a generating set of $I_L$ when $L=L_{pure}$ might have a common monomial factor according to Theorem~\ref{Markov_pure}. Note also that if $E=\{\ub_1,\ldots,\ub_r\}$  is a basis of $L$ such that $\ub_1\in L^+$  and $\supp(\ub_1)=\sigma$ then the set $\{1-x^{\ub_1}, x^{\ub_2^+}-x^{\ub_2^-},\ldots, x^{\ub_r^+}-x^{\ub_r^-}\}$ is a cardinality-minimal generating set of $I_L$. On the other hand such a basis $E$ for $L=L_{pure}$ exists by Corollary~\ref{cor_basis_pure}. This implies that if $L=L_{pure}$, the ideal $I_L$ is always a complete intersection, which follows from \cite[Theorem 2.1]{ES}. In section 5 we determine when $I_L$ is a binomial complete intersection ideal for arbitrary lattices. In the next example we show that all conditions of Theorem~\ref{Markov_pure} are necessary.

\begin{ex1}\label{all_3_needed} 
{\em Let $L\subset \Z^3$ be the lattice generated by $(1,1,0),(0,5,0)$. It is not hard to see that $I_L=\langle 1-x_1x_2,1-x_1^5\rangle\subset\Bbbk[x_1,x_2,x_3]$. Consider also the following sets: $E_{23}=\{x_1^5-1, x_2^5-1\}$, $E_{13}=\{x_1x_2-1, x_2^5x_3-x_3\}$, $E_{12a}=\{x_1^5-1, x_2^2x_1-x_2\}$ and $E_{12b}=\{x_1^2x_2^2-x_1x_2, x_1^5x_2-x_2\}$. One can show  via Gr\" obner bases computation that these sets do not generate $I_L$. We note that each of these sets satisfies all except one of the conditions of Theorem~\ref{Markov_pure}, the missing index indicating which one.} 
\end{ex1}

\subsection{Minimal generating sets of lattice ideals}
Next we characterize the generating sets of lattice ideals starting from the criterion given in \cite[Introduction]{HM}. Let $L\subset \Z^n$ be a lattice and $S$ a subset of $I_L$ consisting of binomials of the form $x^{\ub^+}-x^{\ub^-}$ where $\ub\in L$. Let $F$ be an $I_L$-fiber. The sequence
 $ (x^{\ab_1}, x^{\ab_2},\ldots, x^{\ab_k})$ is an {$S$-\it path}  from $x^{\ub}$ to $x^{\vb}$ if 
\begin{itemize}
\item{} $x^{\ab_1}=x^{\ub}$, $x^{\ab_k}=x^{\vb}$
\item{} for $j=1,\ldots, k$  each  $x^{\ab_j}$ in the sequence belongs to the fiber $F$ and
\item{} $x^{\ab_j}-x^{\ab_{j+1}}$ is equal to $x^{\wb_j} B_j$ or $-x^{\wb_j} B_j$ for some $B_j\in S$, $\wb_j\in \mathbb{N}^n$.
\end{itemize}

\begin{thm1}[\cite{HM}]\label{criterion_generating_set} 
The set $S$ of binomials of $I_L$  is a generating set of $I_L$  if and only if for every $I_L$-fiber $F$ there is an $S$-path between any two elements of $F$.
\end{thm1}

Let $F$ be an $I_L$-fiber and $G(M_F)=\{ x^{\ab_1},\ldots, x^{\ab_s}\}$. We define a relation ``$\sim$" among the elements of $G(M_F)$ as follows:
\[
x^{\ab_i}\sim x^{\ab_j} \textrm{ iff } (\ab_i +L^+)\  \bigcap \ (\ab_j +L^+) \neq  \emptyset. 
\]

We note that if $L^+ =\{\bf 0\}$ then $x^{\ab_i}\sim x^{\ab_j}$  only when $x^{\ab_i}= x^{\ab_j}$.

\begin{lem1} 
``$\sim$" is an equivalence relation among the elements of $G(M_F)$.
\end{lem1}

\begin{proof} It is enough to show transitivity. We can assume that $L_{pure}\neq \{\bf 0\}$, the other case being trivial. Suppose that 
\[
x^{\ab_i}\sim x^{\ab_j} \textrm{ and }  x^{\ab_j}\sim x^{\ab_k} \ .
\]
Thus there exist $\ub_i, \ub_j, \vb_j, \vb_k \in L^+$ such that
\[
\ab_i+\ub_i=\ab_j+\ub_j \textrm{ and } \ab_j+\vb_j=\ab_k+\vb_k\ .
\]
Therefore
\[ 
\ab_i+(\ub_i+\vb_j)=\ab_k+(\ub_j+\vb_k),\ 
\]
and the proof is complete.
\end{proof}

 For the following lemma, we recall that $\ub^{\sigma}$ stands for the vector $(u_i)_{i\notin \sigma}$. 

\begin{lem1}\label{same_vertices_graph_lem}
Let $G(M_F)=\{ x^{\ab_1},\ldots, x^{\ab_s}\}$. The following holds for the elements of $G(M_F)$: 
\[
 x^{\ab_i}\sim x^{\ab_j} \textrm{ if and only if  }\ab_i^\sigma=\ab_j^\sigma \ .
\]
\end{lem1}

\begin{proof} We can assume that $L^+ \neq \{\bf 0\}$, the other case being trivial. Let $\wb\in L^+$ such that $\supp (\wb)=\sigma$. Suppose that $\ab_i^\sigma=\ab_j^\sigma$. Since $\ub=\ab_i-\ab_j\in L$  it follows that  $\ub^\sigma={\bf 0}$ and thus $\supp (\ub)\subset \sigma$. Therefore there exists  $\lambda\gg 0$ such that  $\ub+\lambda \wb\in \N^n$. Since $\ub+\lambda\wb\in L^+$ and $\ab_i+\lambda \wb=\ab_j+(\ub+\lambda \wb)$ it follows that $x^{\ab_i}\sim x^{\ab_j}$.

Suppose now that $x^{\ab_i}\sim x^{\ab_j}$. There exist $\ub_i, \ub_j \in L^+$ such that $\ab_i+\ub_i=\ab_j+\ub_j$. Therefore
\[
\ab_i^\sigma=(\ab_i+\ub_i)^\sigma=(\ab_j+\ub_j)^\sigma=\ab_j^\sigma\ ,
\]
and we are done.
\end{proof}

\begin{lem1}\label{distinct_vertices_graph_lem} Let $G(M_F)=\{ x^{\ab_1},\ldots, x^{\ab_s}\}$. The following holds for the elements of $G(M_F)$:
\[
 x^{\ab_i} \not\sim x^{\ab_j} \textrm{ if and only if  } \ab_i^\sigma , \ab_j^\sigma \textrm{  are incomparable.}
\]
\end{lem1}

\begin{proof} We will show that $x^{\ab_i} \not\sim x^{\ab_j}$ implies that  $\ab_i^\sigma $ and $\ab_j^\sigma$ are incomparable. Suppose otherwise. Without loss of generality we can assume that $\ab_i^\sigma-\ab_j^\sigma>{\bf 0}$. Let $\ub=\ab_i-\ab_j$. Let $\wb\in L^+$ such that $\supp(\wb)=\sigma$. We can find $\lambda\in \N$ large enough so that all coordinates $(u+\lambda w)_i>0$ for $i\in \sigma$.   Since $\lambda \wb^{\sigma}={\bf 0}$ and $\ub^\sigma>{\bf 0}$  it follows that  $\ub+\lambda \wb>{\bf 0}$ and thus $\ub+\lambda \wb\in L^+$. Since $\ab_i+\lambda \wb=\ab_j+(\ub+\lambda \wb)$ and $\lambda \wb$, $\ub+\lambda \wb \in L^+$ it follows that $x^{\ab_i}  \sim x^{\ab_j}$, a contradiction. The other implication follows immediately from the previous lemma.
\end{proof}

Lemmas \ref{same_vertices_graph_lem} and \ref{distinct_vertices_graph_lem} imply that there are only two possibilities for $\ab_i^\sigma, \ab_j^\sigma$, when $x^{\ab_i}, x^{\ab_j}$ are minimal generators of $M_F$: either $\ab_i^\sigma=\ab_j^\sigma$ or $\ab_i^\sigma, \ab_j^\sigma$ are incomparable. If $X\subset \Z^n$  by $X^{\sigma}$ we mean the set whose elements  are the vectors $\ub^\sigma$ where $\ub\in X$. In particular  
\[
 G(M_F)^\sigma=\{x^{\ub^{\sigma}}: x^{\ub}\in G(M_F)\}\ .
\]
Note that the cardinality of $G(M_F)^\sigma$ might be less than the cardinality of $G(M_F)$.

\begin{lem1}\label{vertices_equiv_fibers_lem} 
If $F, F'$ are two equivalent $I_L$ fibers, then $G(M_F)^\sigma = G(M_{F'})^\sigma$.
\end{lem1}

\begin{proof}
Since $F, F'$ are  equivalent $I_L$-fibers, there exist monomials $x^{\ub},x^{\vb}$ such that $x^{\ub}F\subset F'$ and $x^{\vb} F'\subset F$. Therefore   $x^{\ub+\vb}F\subset F$ and  $\ub+\vb\in L^+$.  Since $\ub,\vb \in \N^n $ it follows that $\supp(\ub),\supp(\vb)\subset\sigma$. Let $G(M_F)=\{ x^{\ab_1},\ldots, x^{\ab_s}\}$, $G(M_{F'})=\{ x^{\bb_1},\ldots, x^{\bb_r}\}$. To show the desired equality it suffices to show that for any $i\in[s]$ there is $j\in [r]$ so that  $\ab_i^\sigma=\bb_j^\sigma$, the other inclusion being taken care by symmetry.

Since $x^{\ab_i}x^{\ub}$ is in $F'$ there exists   $j\in [r]$ such that $x^{\bb_j}$ divides $x^{\ab_i}x^{\ub}$. Therefore $\ab_i+\ub-\bb_j\in \N^n$. Since $\supp(\ub)\subset \sigma$ it follows that 
\[
(\ab_i+\ub-\bb_j)^\sigma= \ab_i^\sigma-\bb_j^\sigma \ge {\bf 0} \ 
\]
and $\ab_i^\sigma \ge \bb_j^\sigma$. Similarly there exists   $k\in [s]$ such that $\bb_j^\sigma \ge \ab_k^\sigma$. Therefore $\ab_i^\sigma \ge \bb_j^\sigma \ge \ab_k^\sigma$ and $\ab_i^\sigma \ge \ab_k^\sigma$.  It follows that $\ab_i^\sigma =\ab_k^\sigma$ and thus $\ab_i^\sigma=\bb_j^\sigma$.
\end{proof}

Let $F$ be any $I_L$-fiber. We construct a graph $G_{\overline F}$ and then we build  on $G_{\overline F}$ to construct a graph $\Gamma_{\overline F}$, that will be crucial in determining when a set of binomials of $I_L$ generates $I_{L,{\le \overline F}}$.

\begin{def1}\label{G_graph}
{\rm Let  $F$ be an $I_L$-fiber, $G(M_F)^\sigma=\{ x^{\ab^\sigma_1},\ldots, x^{\ab^\sigma_k}\}$ where $x^{\ab_i} \in G(M_F)$ for $i\in [k]$. We define  $G_{\overline F}=(V(G), E(G))$ to be the graph with $V(G)=[k]$,  and 
\[
E(G)=\{ \{i,j\}:\   \exists\ x^{\ub_i}, x^{\ub_j}\in F \textrm{ such that } \ub_i^\sigma=\ab_i^\sigma,\ \ub_j^\sigma=\ab_j^\sigma, \ x^{\ub_i}-x^{\ub_j}\in I_{L,{< \overline F}}\}\ .
\]}
\end{def1}

The graph $G_{\overline F}$ is independent of the fiber representative $F$, up to reordering of the vertices. This is immediate for  $V(G)$,  by Lemma~\ref{vertices_equiv_fibers_lem}. Next we show the independence of  $E(G)$. Suppose that $x^{\ub_i}-x^{\ub_j}\in I_{L,{< \overline F}}$ where $x^{\ub_i}, x^{\ub_j}\in F$ and $\ub_i^\sigma=\ab_i^\sigma$, $\ub_j^\sigma=\ab_j^\sigma$. Let $F'\in \overline F$ and let $x^{\ub},x^{\vb}$ be such that $x^{\ub} F\subset F'$ and $x^{\vb} F'\subset F$. By Lemma \ref{prop_fiber_mult},  $\supp(\ub)\cup\supp(\vb)\subset \sigma$ and thus $\ub^\sigma=\vb^\sigma={\bf 0}$. Moreover $x^{\ub} x^{\ub_i}-x^{\ub} x^{\ub_j} \in I_{L,{< \overline F}}$ , $({\ub+\ub_i})^{\sigma}= \ub^\sigma+\ub_i^\sigma=\ub_i^\sigma=\ab_i^\sigma$,   $({\ub+\ub_j})^{\sigma}=\ab_j^\sigma$ and thus $E(G)$ is  independent on the choice of the fiber representative $F$.

\begin{ex1}\label{graph_G_F}
{\em Let  $L$ be the sublattice of $\Z^5$, generated by $\vb_1=(3,0,1,-1,0)$, $\vb_2=(0,1,6,0,-1)$, $\vb_3=(1,1,0,0,0)$ and  $\vb_4=(5,0,0,0,0)$. It is easy to see that  $L_{pure}=\Z\vb_3+\Z\vb_4$ and thus $\sigma=\{1,2\}$. Since $\overline{F_{x_5}}=\overline{F_{x_4^6}}$, the $I_L$-fiber  $F=F_{x_5}$ is a Markov fiber\footnote{This is  proved in detail in the last section.}. 
A straightforward computation shows that \[G(M_{F})=\{x_5,x_1^4x_3^6,x_1x_3^5x_4,x_1^3x_3^4x_4^2,x_3^3x_4^3,x_1^2x_3^2x_4^4,x_1^4x_3x_4^5,
x_1x_4^6,x_2x_3^6, x_2^4x_3^5x_4,x_2^2x_3^4x_4^2,\] \[x_2^3x_3^2x_4^4,x_2x_3x_4^5, x_2^4x_4^6\}.\] 
Indeed, if one solves the system that results from the observation that $x_1^{v_1}\cdots x_5^{v_5} \in F$ if and only if 
$$(v_1,v_2,v_3,v_4,v_5-1)=\beta_1\vb_1+\cdots+\beta_4\vb_4=(3\beta_1+\beta_3+5\beta_4,\beta_2+\beta_3,\beta_1+6\beta_2,-\beta_1,-\beta_2),$$
for some $(\beta_1,\ldots,\beta_4)\in\Z^4$, then obtains the above displayed monomials.
Therefore 
\[
G(M_{F})^{\sigma}=\{x_5,x_3^6,x_3^5x_4,x_3^4x_4^2,x_3^3x_4^3,x_3^2x_4^4,x_3x_4^5,x_4^6\}=
\{x_5\}\cup \{x_3^{6-k}x_4^k: \ k=0,\ldots, 6\}.
\]
Thus the graph $G_{\overline {F}}$ consists of $8$ vertices. One can show  that the vertex that corresponds to $x_5$ is isolated.
If $x_5$ was not isolated,  there would be $x^{\ub},x^{\vb}\in F$, such that $x^{{\ub}^{\sigma}}=x_5$, $x^{{\vb}^{\sigma}}=x_3^{6-k}x_4^k$ and
 $x^{\ub}-x^{\vb}\in I_{<\overline{F}}$.  But this means that there would be  a $\wb\in L^+$,  so that $x^{\ub}=x^{\wb}x_5$. Thus, since  $x_5(x^{\wb}-1) \in  I_{<\overline{F}}$,  we get that $x_5-x^{\vb}=(x^{\ub}-x^{\vb})-x_5(x^{\wb}-1)  \in I_{<\overline{F}}$.
This leads to a contradiction,  by Theorem~\ref{weak_betti}, since $F_{x_5}$ is a Markov fiber. It is also relatively easy to show that any other two vertices of  $G_{\overline {F}}$  are connected by an edge. Consider for example the vertices that correspond to the monomials $x_3^6$ and $x_3^5x_4$. They are connected via an edge since $(x_1^4x_3^6)^\sigma=x_3^6$, $(x_1x_3^5x_4)^\sigma=x_3^5x_4$,  $x_1^4x_3^6-x_1x_3^5x_4=x_1x_3^5(x_1^3x_3-x_4)$ and $ \overline F_{x_4}<_{I_L} \overline F_{x_5}$, implying  that  $x_1^4x_3^6-x_1x_3^5x_4 \in I_{<\overline{F}}$.
Therefore $G_{\overline {F}}$ has two connected components: an isolated vertex corresponding to $x_5$ and the complete graph on the remaining seven vertices.}
\end{ex1}

\medskip
\begin{def1}\label{gamma_graph}
{\rm We let $\Gamma_{\overline F}$ to be  the complete graph  whose vertices are the connected components of $G_{\overline F}$. Let $B=x^{\ub}-x^{\vb}\in I_L$ such that $F_B\in \overline F$. We identify $B$ with an edge of $\Gamma_{\overline F}$ if $x^{\ub^\sigma}\neq x^{\vb^\sigma}$, $B\in I_{L,{\le \overline F }}$ and  $B\not\in   I_{L,{<\overline F}}$. For a subset $S$ of binomials of $I_L$ we denote by $\Gamma_{\overline F}(S)$ the subgraph of $\Gamma_{\overline F}$ induced by the binomials $B\in  S$ such that $F_B\in \overline F$.}
\end{def1}

We note that  different binomials might correspond to the same edge of  $\Gamma_{\overline F}$.

\begin{lem1}\label{lem_main_thm}
Let $L\subset \Z^n$ be a lattice and $S$ a subset of $I_L$ consisting of binomials such that  $I_{L_{pure}}=\langle S_{\overline F_{\{1\}}}\rangle$ and $\Gamma_{\overline F}(S)$ is  a spanning tree of $\Gamma_{\overline F}$ for every $I_L$-fiber $F$. Then the set  $S$ is a generating set of $I_L$. 
\end{lem1}

\begin{proof}
If $L=L_{pure}$ then the conclusion is straightforward since $G_{\overline F_{\{1\}}}$ is an isolated vertex, and so is $\Gamma_{\overline F_{\{1\}}}$. Thus we may assume that the lattice $L$ is non-pure or equivalently that there exists a fiber $F$ such that $F\notin\overline F_{\{1\}}$. We will show that for any $I_L$-fiber $F$ and any $x^{\ub}$, $x^{\vb}\in F$ there is an $S$-path between $x^{\ub}$, $x^{\vb}$. This was already noticed if $F\in \overline F_{\{1\}}$. By Theorem \ref{chain_least} we can assume that there is an $S$-path between any two elements of $G$ for all $G$ such that $\overline G<_{_{I_L}} \overline F$.  We note that  $\ub-\vb\in L$. Suppose that $G(M_F)^\sigma=\{ x^{\ab^\sigma_1},\ldots, x^{\ab^\sigma_k}\}$ where $x^{\ab_1},\ldots, x^{\ab_k}\in G(M_F)$. We examine three cases.

{\bf Case 1.} If $\ub^\sigma=\vb^\sigma$ then  $(\ub-\vb)^\sigma={\bf 0}$ and by Corollary \ref{pure_supp_outside} it follows that  $\ub-\vb\in L_{pure}$. Therefore $x^{\ub}-x^{\vb}\in I_{L_{pure}}$ and  since  $I_{L_{pure}}=\langle S_{\overline F_{\{1\}}}\rangle$ it follows that $x^{\ub}-x^{\vb}\in  \langle S\rangle$.

{\bf Case 2.} Suppose that  $\ub^\sigma\neq \vb^\sigma$ and that the vertices of $G_{\overline{F}}$ corresponding to $\ub^\sigma$ and $\vb^\sigma$ are in the same connected component  of $G_{\overline{F}}$. Assume that $\ub^\sigma=\ab_i^\sigma$ and $\vb^\sigma=\ab_j^\sigma$ and that $i=i_1,\ldots,i_l=j$ is a path in $G_{\overline{F}}$. By applying induction on $l$ it is enough to prove the statement when $l=2$ and $\{i,j\}$ is an edge of $G_{\overline{F}}$. It follows from the definition of $G_{\overline{F}}$ that there exists a binomial $x^{\wb}-x^{\zb}\in I_{L,{<\overline F}}$ such that $\wb^\sigma=\ub^\sigma$, $\zb^\sigma=\vb^\sigma$. Moreover $x^{\wb}, x^{\zb}\in G$ where $\overline G<_{_{I_L}} \overline F$. Thus there is a monomial $x^{\ab}$ such that $x^{\ab}G\subset F$. By Case 1 above, there is an  $S$-path from $x^{\ub}$ to $x^{\wb+\ab}$ and an $S$-path from $x^{\vb}$ to $x^{\zb+\ab}$. By assumption there is an $S$-path from $x^{\wb}$ to $x^{\zb}$, and thus also from $x^{\wb+\ab}$ to $x^{\zb+\ab}$. Putting the $S$-paths together one gets  an $S$-path from $x^{\ub}$ to $x^{\vb}$. We point out that the above argument shows  that $x^{\ub}-x^{\vb}\in I_{L,{<\overline F}}$.

{\bf Case 3.} Suppose that  $\ub^\sigma\neq \vb^\sigma$ and that the vertices of $G_{\overline{F}}$  corresponding to $\ub^\sigma$ and $\vb^\sigma$ are in disconnected components of $G_{\overline{F}}$. Since $S$ determines a spanning tree of $\Gamma_{\overline F}$ there is a series of edges in $\Gamma_{\overline F}$ that leads from the component that corresponds to $x^{\ub^\sigma}$ to the component that corresponds to $x^{\vb^\sigma}$. As before it is enough to prove the statement when the components are adjacent. This means that there exists a binomial $B= x^{\ub'}-x^{\vb'}\in (I_{L,{\le \overline F}}\setminus I_{L,{<\overline F}})\cap S$ such that  ${\ub'}^\sigma$, $\ub^\sigma$ correspond to the same connected component of $G_{\overline{F}}$ and similarly for  ${\vb'}^\sigma$, $\vb^\sigma$. The monomials $x^{\ub'}$,$x^{\vb'}$ of $B$, belong to a fiber  equivalent to $F$. It follows  that there is $\bb\in \N^n$ such that $x^{\ub'+\bb}$ and $x^{\vb'+\bb}$ belong to $F$ and thus the sequence  $(x^{\ub'+\bb},x^{\vb'+\bb})$ is  an $S$-path from $x^{\ub'+\bb}$ to $x^{\vb'+\bb}$. By Case 2 above, there is an $S$-path from  $x^{\ub}$ to $x^{\ub'+\bb}$ and an $S$-path from $x^{\vb}$ to $x^{\vb'+\bb}$. Joining these paths one gets an $S$-path from $x^{\ub}$ to $x^{\vb}$.
\end{proof}

\noindent We let $t(\overline F)$ denote the number of vertices of $\Gamma_{\overline F}$. Thus
\[
 t(\overline F):= | V(\Gamma_{\overline F})|. 
\]
We note that to construct a spanning tree of $\Gamma_{\overline F}$ we need exactly $t(\overline F)-1$ binomials. To prove the next theorem we will use Theorem \ref{criterion_generating_set}.

\begin{thm1}\label{main_thm} 
Let $L\subset \Z^n$ be a lattice and $S$ a subset of $I_L$ consisting of binomials. The set $S$ is a cardinality-minimal generating set of $I_L$ if and only if the following conditions are satisfied:
\begin{itemize}
\item{} $\Gamma_{\overline F}(S)$ is a spanning tree of $\Gamma_{\overline F}$ for every $I_L$-fiber $F$, and $|S_{\overline F}|=t(\overline F)-1$ for every $I_L$-fiber $F$ such that $F\notin \overline F_{\{1\}}$,
\item{} $|S_{\overline F_{\{1\}}}|= \rank (L_{pure})$ and
\item{} $I_{L_{pure}}=\langle S_{\overline F_{\{1\}}}\rangle$.
\end{itemize}
\end{thm1}

\begin{proof}
Suppose  that $S$ is a cardinality-minimal generating set of $I_L$. We will show that $S$ satisfies the three conditions of the theorem. We note that since $S$ is a generating set of $I_L$, by Proposition~\ref{leq_prop} and Remark~\ref{fiber_of_one} it follows that $\langle S_{\overline F_{\{1\}}}\rangle =I_{L_{pure}}$. Moreover if  $|S_{\overline F_{\{1\}}}|> \rank (L_{pure})$, then by Theorem~\ref{Markov_pure} one can replace the binomials in $S_{\overline F_{\{1\}}}$ by a cardinality-minimal generating set of $I_{L_{pure}}$. The new set thus produced is  still a generating set of $I_L$, according to Lemma~\ref{lem_main_thm}, and has smaller cardinality than $S$, a contradiction. Next we show that for an arbitrary $I_L$-fiber $F$,  $\Gamma_{\overline F}(S)$ is a spanning tree of $\Gamma_{\overline F}$. Indeed, by Theorem~\ref{criterion_generating_set}, $S$ induces a spanning graph in $F$. Since $G_{\overline F}$ comes from $F$ by identifying components and similarly for $\Gamma_{\overline F}$  from $G_{\overline F}$, it follows that $\Gamma_{\overline F}(S)$ is  a spanning graph of $\Gamma_{\overline F}$. We show that $\Gamma_{\overline F}(S)$ is a  tree of $\Gamma_{\overline F}$. If not, $\Gamma_{\overline F}(S)$ has a cycle in $\Gamma_{\overline F}$. We omit from $S$ the binomial that induces an edge on this cycle. The resulting set still satisfies the conditions of  Lemma~\ref{lem_main_thm} and is thus a generating set of $I_L$ of smaller cardinality, a contradiction. Similarly, if for some fiber $F$ such that $F\notin \overline F_{\{1\}}$ we have $|S_{\overline F}|> t(\overline F)-1$ then there is a binomial in $S_{\overline F}$ that does not correspond to an edge of $\Gamma_{\overline F}(S)$ or two binomials that correspond to the same edge. Then one binomial could be omitted from $S$ and the resulting set would still be a generating set of $I_L$.

Conversely let $S$ be a set that satisfies the three conditions of the theorem. By Lemma~\ref{lem_main_thm}, $S$ is a generating set of $I_L$. Suppose that 
there is a cardinality-minimal generating set $S'$ of $I_L$ such that $|S'|<|S|$. Let $F$ be a  Markov fiber  such that $|S'_{\overline F}|< |S_{\overline F}|$.  We note that $F\notin  \overline F_{\{1\}}$ since $|S'_{\overline F_{\{1\}}}|\ge \rank (L_{pure})=|S_{\overline F_{\{1\}}}|$. Moreover $|S_{\overline F}|=t(\overline F)-1\leq  |S'_{\overline F}|$, a contradiction. Therefore $S$ is a cardinality-minimal generating set of $I_L$. 
\end{proof}

\begin{rem1}\label{positive_main_thm}
{\rm When $L_{pure}=\{\bf 0\}$ then the construction and conditions on $\Gamma_{\overline F}$ coincide with the construction and conditions of $\mathcal{S}_{\bf b}$ from \cite{CKT}, since by Proposition~\ref{infinitefiber} $G(M_F)=F$. Furthermore, in this case Theorem~\ref{main_thm} becomes the join of \cite[Theorem 2.6]{CKT} and \cite[Theorem 2.7]{CKT}. To understand better how Theorem~\ref{main_thm} works to compute all (infinitely many) cardinality-minimal generating sets in the general case of a non-positive lattice, see the example worked in great detail from Section 6.}
\end{rem1}

\noindent We remark that for all but finitely many equivalence classes of fibers $\overline F$, $t(\overline F) =1$. Indeed by Corollary~\ref{markov_fib_inv_cor} the set consisting of equivalence classes of Markov fibers is finite. If an $I_L $-fiber $F$   is not a Markov fiber, then by Theorem~\ref{weak_betti}
it follows that $I_{L,{\leq\overline F}}=I_{L,{< \overline F}}$ and hence $G_{\overline F}$ consists of only one connected component. The next result is the main theorem of this section. Its proof is an immediate consequence of Theorem~\ref{main_thm}. 

\begin{thm1} \label{cor_rank+graph}
Let $L\subset \Z^n$ be a lattice and $\mu(I_L)$ be the cardinal of a cardinality-minimal generating set of $I_L$. Then 
\[
 \mu(I_L)= \rank (L_{pure})+\sum_{ \overline F\neq \overline{F}_{\{1\}}} (t(\overline F)-1),\ 
\]
where the sum runs over all distinct equivalence classes of Markov fibers.
\end{thm1}

\noindent  The next corollary follows from Corollary  \ref{markov_fib_inv_cor} and Theorem  \ref{main_thm}. It generalizes the corresponding result for positive lattices, see \cite[Theorem 2.5]{BCMP} and \cite[Theorem 1.3.2]{DSS}.

\begin{cor1}\label{multiset_inv} 
Let $L\subset \Z^n$ be a lattice, $\mu=\mu(I_L)$ and  $\{B_1,\ldots, B_{\mu}\}$  a cardinality-minimal generating set of $I_L$. The multiset
\[
\{ {\overline F_{B_1}},\ldots, {\overline F_{B_{\mu}}} \}
\]
is an invariant of $I_L$.
\end{cor1}

\noindent The following result concerns an arbitrary minimal generating set of $I_L$: the play in the cardinality of such a set only concerns the pure part of $L$. Note that if $L_{pure}\neq 0$ then $I_L$ can be minimally generated by $\mu(I_L)+k$ binomials for every $k\in\N$. More precisely, all extra $k$ binomials are given by the fiber $F_{\{1\}}$ with a similar argument as the one given in introduction for the lattice $L=\langle (1,1), (5,0)\rangle$. The proof of Corollary \ref{cor_multiset} follows directly from the proof of Theorem  \ref{main_thm}.

 \begin{cor1}\label{cor_multiset} Let $L$ be a lattice and $S=\{B_1,\ldots, B_t\}$ a minimal generating set of $I_L$. The multiset
\[ 
\{ {\overline F_{B_i}}:\  {{F_{B_i}}} \not\in \overline F_{\{1\}}\}
\]
 is an invariant of $I_L$.
 \end{cor1}

\subsection{Applications}
In the last years, due to applications of minimal generating sets to algebraic statistics, there is an interest in determining the {\it indispensable binomials} of a lattice ideal, see \cite{Hi-O, HO, CKT, ATY, OVT}. Of particular interest is the case when all elements in a  minimal generating set of the lattice ideal are indispensable as is the case for generic lattice ideals, see \cite{PS}.

\begin{def1}\label{indisp_weak_indisp}
{\rm A binomial is called {\it  indispensable} ({\it weakly indispensable}, respectively) if it appears in every minimal generating set (every cardinality-minimal generating set, respectively) of $I_L$ up to a constant multiple. A monomial $x^{\ub}$ is called {\it indispensable} ({\it weakly indispensable}, respectively) if for every minimal generating set (every cardinality-minimal generating set, respectively) $S$ of $I_L$ there is a binomial $B\in S$ so that $x^{\ub}$ is a monomial term of $B$.}
\end{def1}

Note that an indispensable binomial/monomial is weakly indispensable. If $L$ is positive then the characterization of indispensable binomials is given by \cite[Corollary 2.10]{CKT}, see also \cite[Theorem 3.4]{CKT}. Note that in the case of positive lattices the notions of indispensable and weakly indispensable coincide. Recently, a polynomial-time algorithm for computing all indispensable binomials from a given system of binomial generators of an arbitrary binomial ideal was given in \cite[Algorithm 1]{CTVJSC}. Next we complete the classification of indispensable binomials/monomials of lattice ideals in the case of non-positive lattices. 


\begin{thm1}\label{thm_indisp_bin} 
Let $L$ be a non-positive lattice. Then $I_L$ has no indispensable binomials, but it has one indispensable monomial, $1=x^{\bf 0}$. If $\rank (L_{pure})>1$ then there are no weakly indispensable binomials and only one weakly indispensable monomial $1=x^{\bf 0}$. If $\rank (L_{pure})=1$ there exists exactly one weakly indispensable binomial and exactly two weakly indispensable monomials. 
\end{thm1}
\begin{proof}
Let $L$ be a lattice such that $\rank (L_{pure})=r\geq 1$ and $S$ a cardinality-minimal generating set of $I_L$. In order to determine the (weakly) indispensable binomials/monomials of $I_L$ we analyze the two cases of Markov fibers which might contain them: 1) $F\in\overline F_{\{1\}}$ and 2) $F\notin\overline F_{\{1\}}$. 

In the first case,  if $L_{pure}$ has rank 1 then $L_{pure}=\langle \ub\rangle$, where $\ub\in L^+$ is $L$-primitive and by Theorem~\ref{Markov_pure} and Theorem~\ref{main_thm} we have that  $S_{\overline F_{\{1\}}}=\{x^{\ub}-1\}$. Thus, when $\rank (L_{pure})=1$ we have one weakly indispensable binomial and two weakly indispensable monomials, which correspond to fiber $F_{\{1\}}$. Replacing $x^{\ub}-1$ from a cardinality-minimal generating set with $x^{2\ub}-1,x^{3\ub}-1$ we obtain a  minimal generating set of $I_L$. Hence we have no indispensable binomials corresponding to a fiber $F\in\overline F_{\{1\}}$. However, $1$ is an indispensable monomial of $I_L$ because otherwise the ideal $I_L$ would be contained in the maximal ideal $\langle x_1,\ldots,x_n\rangle$, a contradiction since $\rank (L_{pure})>0$. 

Suppose now that  $L_{pure}$ has rank $r>1$. Without loss of generality we can assume that $S_{\overline F_{\{1\}}}=\{x^{\ub_1}-1,\ldots,x^{\ub_r}-x^{\vb_r}\}$ satisfies the three conditions of Theorem~\ref{Markov_pure}. For every $i\geq 2$ let $\ub'_i=\ub_i+\ub_1$ and $\vb'_i=\vb_i+\ub_1$ and note that $\ub'_i-\vb'_i=\ub_i-\vb_i$. By Theorem~\ref{Markov_pure} it follows  that $\{x^{\ub_1}-1,x^{\ub'_2}-x^{\vb'_2},\ldots,x^{\ub'_r}-x^{\vb'_r}\}$ is also a cardinality-minimal generating set of $I_{L_{pure}}$ having only one binomial in common with $S_{\overline F_{\{1\}}}$. Since $\rank(L_{pure})>1$ there are infinitely many $L$-primitive elements of full support and thus infinitely many bases of $L_{pure}$ as in Corollary~\ref{cor_basis_pure}. Any of these bases induces a cardinality-minimal generating set of $I_{L_{pure}}$, see Theorem~\ref{Markov_pure}. By applying the above trick to any of these cardinality-minimal generating sets we conclude that there are no weakly indispensable binomials for $I_L$ corresponding to a fiber $F\in\overline F_{\{1\}}$, but there is one weakly indispensable monomial of $I_L$ belonging to $F_{\{1\}}$ and is $1=x^{\bf 0}$. Obviously there are no indispensable binomials in this case, but $1$ is an indispensable monomial with the same argument as before. 

In the second case assume that there exists a Markov fiber $F$ such that $F\notin \overline F_{\{1\}}$, that is $\overline F>_{_{I_L}}\overline F_{\{1\}}$. Then we prove that there are  infinitely many distinct choices for the binomials that determine any edge in a spanning tree of $\Gamma_{\overline F}$. Since $F$ is a Markov fiber,  Theorem~\ref{weak_betti} says that  $I_{L,{<\overline F}}\neq I_{L,{\leq \overline F}}$. Thus there exists   $B=x^{\ub}-x^{\vb}$, such that $F_B\in \overline F$ and $B\notin I_{<\overline F}$. We note that ${\ub}^{\sigma}\neq {\vb}^{\sigma}$:  otherwise $\ub-\vb\in L_{pure}$ and $x^{\ub}-x^{\vb}\in I_{L_{pure}}=I_{L,{\leq\overline F_{\{1\}}}}\subset I_{L,{<\overline F}}$, a contradiction. Therefore $B$ produces  an edge in $\Gamma_{\overline F}$ and can be made part of a cardinality-minimal generating set $S$ of $I_L$. Since $L$ is non-positive there exist vectors $\wb_1,\wb_2\in L^+\setminus\{{\bf 0}\}$. Then $B'=x^{\ub+\wb_1}-x^{\vb+\wb_2}$ gives exactly the same edge as $B$ and can replace $B$ in $S$. Hence we obtain that there are no weakly indispensable binomials or monomials of $I_L$ corresponding to a fiber $F\notin\overline F_{\{1\}}$. This also implies that in the second case there are no indispensable binomials or monomials of $I_L$.     
\end{proof}

One nice practical application of Theorem~\ref{thm_indisp_bin} is to detect whether a lattice is positive or not. Indeed, given any arbitrary set of vectors which span a lattice $L$ then: 1) apply one of the existing algorithms, for example the one from \cite{BSR} or \cite{HM}, to compute a generating set of $I_L$, and 2) decide that $L$ is positive if and only if $1$ is not in the support of any binomial from the computed generating set. The latter is true by Theorem~\ref{thm_indisp_bin} since for non-positive lattices the corresponding lattice ideals have $1$ as indispensable monomial while for positive lattices the fiber of $1$ consists only of $1$. Next we present another application following from the proof of Theorem~\ref{thm_indisp_bin}. First we define the {\it universal Markov basis} of $I_L$ to be the union of all cardinality-minimal generating sets of $I_L$, up to a sign.  

\begin{thm1}\label{graver}
Let $L$ be a lattice. If $\rank (L_{pure})>1$ or $\rank(L_{pure})=1$ and $L\neq L_{pure}$ then the universal Markov basis of $I_L$ is infinite. 
\end{thm1}

Note that the alternative definition of universal Markov basis as the union of all minimal generating sets would also imply the conclusion of Theorem~\ref{graver} and only under the assumption that the lattice is non-positive, as explained in the Introduction. However, both definitions are identical in the case of positive lattices, since a minimal generating set is automatically a cardinality-minimal generating set, and they match the definition of universal Markov basis introduced in \cite[Definition 3.1]{HS}. For positive lattices the vectors corresponding to the universal Markov basis have also an useful algebraic description, see \cite[Proposition 1.4]{CTVJA}. We note that for non-positive lattices $L$ satisfying the hypotheses of Theorem~\ref{graver} the universal Markov basis of $I_L$ is not contained in the Graver basis, since the Graver basis is finite (see \cite[Algorithm 7.2]{St}). Moreover, Theorem~\ref{graver} was the the starting point for classifying the lattices for which the universal Markov basis is contained either in the Graver basis or the universal Gr\" obner basis, see \cite[Theorem 2.3]{CTV} or \cite[Theorem 2.7]{CTV}. 

The Markov complexity of an integer matrix, see \cite{SS} for the definition, is another topic of interest in algebraic statistics, related to the universal Markov basis. We note that the results of this paper and in particular Theorems~\ref{Markov_pure}, \ref{main_thm} are used in \cite[Section 3]{CTV} to give necessary and sufficient conditions for the Markov complexity to be finite.

  
\section{Binomial Complete Intersection Lattice Ideals}

In this section we characterize all binomial complete intersection lattice ideals. This is a problem that engaged mathematicians starting in 1970, see \cite{H}. We recall that a lattice ideal $I_L$ of height $r$ is  a {\it complete intersection} if there exist polynomials $P_1,\ldots, P_r$ such that   $I_L=\langle P_1,\ldots, P_r\rangle$ and $I_L$ is a {\it binomial complete intersection} if there exist binomials $B_1,\ldots, B_r$ such that $I_L=\langle B_1,\ldots, B_r\rangle$. In the case of positive lattices Nakayama's lemma applies and then complete intersection lattice ideals are automatically binomial complete intersections. Furthermore, for positive lattices the problem was gradually solved in a series of papers \cite{H, D,Sta,I,W,N,Sc,RGS,FMS,SSS,MT}. The final conclusion is that $I_L$ is a complete intersection if and only if the   matrix $M$ whose  rows correspond to a basis of $L$ is ``mixed dominating'': every row of $M$ has a positive and negative entry and $M$ contains no square submatrix with this property, see \cite[Theorem 3.9]{MT}. However, when $L$ is non-positive the problem was open and we give a complete answer in Theorem~\ref{compint} and Corollary~\ref{char_binom_compl_inter_cor}.

We recall that $\sigma$ is the maximum support of an element of $L\cap \N^n$, $L_{pure}$ is the sublattice of $L$ generated by the elements of $L\cap \N^n$, $\ub^{\sigma}=(u_i)_{i\notin\sigma}\in(\Z^n)^{\sigma}$ and $L^{\sigma}$ is the sublattice of $(\Z^n)^\sigma$ generated by the vectors $\ub^\sigma$, where $\ub\in L$.

As an immediate consequence of definitions and Corollary~\ref{pure_supp_outside} we have the following equivalences: 
\[
\sigma=\emptyset \ \Leftrightarrow \ L_{pure}=\{{\bf 0}\} \ \Leftrightarrow \ L \text{ is positive } \ \Leftrightarrow \ L^{\sigma}=L,
\]
and
\[
L \text{ is pure } \ \Leftrightarrow \ L=L_{pure} \  \Leftrightarrow \ L^{\sigma}=\{{\bf 0}\}. 
\] 
 First we have the following:

\begin{rem1}\label{pos_graded_lattice}
{\rm  $L^\sigma$ is a positive lattice.}
\end{rem1}
\begin{proof} The conclusion follows at once when $\sigma=\emptyset$. Assume now that $\sigma\neq\emptyset$ and let $\wb\in L\cap \N^n$ be such that $\supp(\wb)=\sigma$. If ${\bf 0}\neq \ub^\sigma\in L^\sigma \cap (\N^n)^\sigma$ then there exists  $k\in \N$, $k\gg 0$ such that $\ub'=\ub+k\wb\in L\cap \N^n$. Thus $\sigma \subsetneq \supp(\ub')$, a contradiction.
\end{proof}

Next we show that the rank of $L$ is determined by the ranks of its sublattice $L_{pure}$ and the lattice $L^\sigma$.

\begin{prop1} \label{lemma}
Let $L\subset \Z^n$ be a lattice. Then $L\cong L^{\sigma}\oplus L_{pure}$.
\end{prop1}
\begin{proof}
The following is a split exact sequence of $\Z$-modules:
\[
0\longrightarrow L_{pure}\hookrightarrow L \stackrel{\pi}{\rightarrow} L^{\sigma} \longrightarrow 0,
\] 
where $\pi$ is the restriction of the canonical projection $\Z^n\rightarrow(\Z^n)^{\sigma}$ given by $\ub\rightarrow \pi(\ub)=(u_i)_{i\notin\sigma}$. For the exactness, we note that $\pi$ is a surjective homomorphism of $\Z$-modules whose kernel is $L_{pure}$, by Corollary~\ref{pure_supp_outside}.
\end{proof}

We consider the lattice ideal $I_{L^\sigma}$ in $R^\sigma:=\Bbbk[x_i: i\not \in \sigma]$. We show that in order to compute the $I_{L^\sigma}$-fibers, it is enough to consider the generating sets of the corresponding $I_L$-fibers.

\begin{lem1}\label{fiber_outside_sigma} 
Let $L\subset\Z^n$ be a lattice and $\ub\in \N^n$. Then the $I_{L^{\sigma}}$-fiber of $x^{\ub^{\sigma}}$ is $G(M_{F_{\ub}})^{\sigma}$, where $F_{\ub}$ is the $I_L$-fiber of $x^{\ub}$.
\end{lem1}
\begin{proof} Let $\ub'=\ub^{\sigma}\in(\N^n)^{\sigma}$ and denote by $F'$ the $I_{L^{\sigma}}$-fiber of $x^{\ub'}$. It follows from Remark \ref{pos_graded_lattice} that $F'$ is finite. We will show that $F'=F_{\ub}^\sigma$ and thus by Proposition~\ref{descfiber} we obtain $F'=G(M_{F_{\ub}})^{\sigma}$, since for any vector ${\bf t}\in L^+$ we have ${\bf t}^{\sigma}={\bf 0}$. To prove this, consider first an element $x^{\vb}\in F_{\ub}$. Then $\vb-\ub\in L$ and $\vb^{\sigma}-\ub^{\sigma}\in L^{\sigma}$. Hence $x^{\vb^{\sigma}}\in F'$ and we obtain the inclusion $F'\supset F_{\ub}^\sigma$. For the converse inclusion, let $x^{\vb'}\in F'$. Then $\ub'-\vb'\in L^{\sigma}$ and it follows that there exists a vector $\wb\in L$ such that $\ub'-\vb'=\wb^{\sigma}$. Therefore $\vb'=(\ub-\wb)^{\sigma}$. Since $\ub-(\ub-\wb)=\wb\in L$ we obtain that $x^{\ub-\wb}\in F_{\ub}$ and $x^{\vb'}\in F_{\ub}^{\sigma}$, as desired.
\end{proof}

Next we show that the cardinality $\mu(I_L)$ of a 
cardinality-minimal generating set of $I_L$ depends on the cardinality 
$\mu(I_{L^\sigma})$ of a cardinality-minimal generating set of $I_{L^\sigma}$.

\begin{thm1} \label{mu} 
Let $L$ be a lattice. Then
\[ 
\mu(I_L)=\mu(I_{L^\sigma})+\rank (L_{pure}).
\]
\end{thm1}
\begin{proof} If $L$ is positive then $L_{pure}=\{{\bf 0}\}$, $L=L^{\sigma}$ and the conclusion is obvious. On the other hand, when $L$ is a pure lattice then $L^{\sigma}=\{{\bf 0}\}$ and the conclusion follows from Theorem~\ref{Markov_pure}. Assume now that $L$ is non-positive and non-pure, that is both $L_{pure}$ and $L^{\sigma}$ are nonzero. By Theorem~\ref{cor_rank+graph} we have $\mu(I_L)=\rank(L_{pure})+\sum_{\overline{F}_{\ub}\neq\overline{F}_{\{1\}}}(t(\overline{F}_{\ub})-1)$. Applying Lemma~\ref{fiber_outside_sigma} the graphs $\Gamma_{\overline F_{\ub}}$ and  $\Gamma_{\overline F_{{\ub}^{\sigma}}}$ are equal for any fiber $F_{\ub}\notin\overline{F}_{\{1\}}$. Thus, applying once more Theorem~\ref{cor_rank+graph} we obtain $\mu(I_{L^{\sigma}})=\sum_{\overline{F}_{\ub}\neq\overline{F}_{\{1\}}}(t(\overline{F}_{\ub})-1)$, since $L^{\sigma}$ is positive. 
\end{proof}

The following theorem follows directly from Proposition~\ref{lemma} and Theorem~\ref{mu} and determines the binomial complete intersection lattice ideals. It is the main theorem of this section.

\begin{thm1} \label{compint}
 Let $L\subset\Z^n$ be a lattice. The ideal $I_L$ is binomial complete intersection if and only if $I_{L^{\sigma}}$ is complete intersection.
\end{thm1}

We can describe the lattices for which $I_{L}$ is a binomial complete intersection. Recall that a mixed dominating matrix $M$ has the property that every row of $M$ has a positive and negative entry and $M$ contains no square submatrix with this property.

\begin{cor1}\label{char_binom_compl_inter_cor}
Let $L\subset\Z^n$ be a lattice. The ideal $I_L$ is binomial complete intersection if and only if either $L^{\sigma}=\{{\bf 0}\}$ or there exists a basis of $L^\sigma$ such that its vectors give the rows of a mixed dominating matrix.
\end{cor1}
\begin{proof} 
If the lattice $L$ is pure, or equivalently $L^{\sigma}=\{{\bf 0}\}$, then by \cite[Theorem 2.1]{ES} or Theorem~\ref{Markov_pure} the ideal $I_L$ is binomial complete intersection. Assume now that $L$ is non-pure. Since by Remark~\ref{pos_graded_lattice} $L^{\sigma}$ is positive, the proof follows from Theorem~\ref{compint} and \cite[Theorem 3.9]{MT}.
\end{proof}

\begin{rem1}
{\rm Let $L\subset \Z^n$ be a non-pure lattice, $r=\rank(L)$ and $r_+=\rank (L_{pure})$. The previous corollary states that $I_L$ is binomial complete intersection if and only if there is a basis of $L$ whose vectors give the rows of
\[
\begin{bmatrix} A&M\cr C&{\bf 0}\end{bmatrix} 
\]
where $A\in \mathcal{M}_{(r-r_+)\times |\sigma|}(\Z)$, the matrix $M\in \mathcal{M}_{(r-r_+)\times (n-|\sigma|)} (\Z)$ is mixed dominating, the matrix $C\in \mathcal{M}_{r_+\times |\sigma|}(\Z)$ is a matrix whose rows satisfy the conditions of Theorem~\ref{Markov_pure}, and ${\bf 0}$ is the zero matrix. For example, any matrix with entries in $\N$ and linearly independent rows has the desired property for $C$ above, while for a mixed dominating matrix  $M$ one can use \cite[Remark 3.17]{MT} or \cite[Theorem 2.2]{FMS}. By working with the appropriate size matrices  one can easily obtain a class of lattice ideals that are binomial complete intersections. }
\end{rem1}

\section{Example}

Algorithms for computing  a generating set for  lattice ideals  were given in \cite{HSt,BSR,HM}. In the following example we use the techniques introduced in this paper to show how to obtain all cardinality-minimal generating sets of a lattice ideal. Let $L$ be the lattice of Example~\ref{graph_G_F}. $L$ is generated by the vectors $\vb_1=(3,0,1,-1,0)$, $\vb_2=(0,1,6,0,-1)$, $\vb_3=(1,1,0,0,0)$ and $\vb_4=(5,0,0,0,0)$. We note that $\vb_5=(0,5,0,0,0)\in L$, 
\[
L^+=L\cap\N^5= \N \vb_3+\N \vb_4+\N \vb_5
\] 
and
\[ 
L_{pure}=  \Z \vb_3+\Z \vb_4 \ .
\]
Thus $\rank (L_{pure})=2$, $\sigma=\{1,2\}$. Let $F$ be an $I_L$-fiber. Since $|\Z^2/ (L_{pure})_\sigma|=5$ it follows by Proposition~\ref{same_nr_equiv} that $\overline F$ consists of five equivalent fibers. In particular, let $F_{\{1\}}$ be the fiber that contains $1=x^{\bf 0}$. As in Example~\ref{example_multi}(b) we see that 
\[
\overline F_{\{1\}} =\{F_{\{1\}} ,F_{x_1},F_{x_1^2},F_{x_1^3},F_{x_1^4}\}\ 
\] 
where for each $0\le k\le 4$, the fiber $F_{x_1^k}$ consists of the monomials $x_1^ix_2^j$ with $i-j\equiv k \mod 5$. Thus if $x^{\ub}\in R$ then 
\[
\overline F_{x^{\ub}}=\{F_{x^{\ub}},F_{x^{\ub}x_1},F_{x^{\ub}x_1^2},F_{x^{\ub}x_1^3},F_{x^{\ub}x_1^4}\}.
\]
Note that every equivalence class of fibers is in the form $\overline F_{x_4^n}$ for some $n \in \N$ and 
\begin{eqnarray} \label{explicit_sequence}
\overline F_{\{1\}}<_{_{I_L}} \overline F_{x_4}<_{_{I_L}} \overline F_{x_4^2}<_{_{I_L}} \cdots <_{_{I_L}} \overline F_{x_4^n}<_{_{I_L}} \cdots .
\end{eqnarray}
 
It is clear that the above equivalence classes of fibers are pairwise distinct: use Lemma~\ref{prop_fiber_mult} and notice that there is no $\wb\in \N^5$ with $\supp(\wb)\subset \sigma$ such that $\wb+(0,0,0,k,0)\in L$. 
 
To find a generating set of $I_L$ in $R=\Bbbk[x_1,\ldots, x_5]$ one first computes the ideal $$ \langle x_4-x_1^3x_3,x_5-x_2x_3^6,1-x_1x_2,1-x_1^5\rangle :(x_1\cdots x_5)^{\infty}$$
using for example CoCoA \cite{Co}. Computing $I_{L,{<\overline  F_{x_4^i}}}$ and $I_{L,{\leq\overline F_{x_4^i}}}$ with the use of Proposition~\ref{leq_prop} we conclude via Theorem~\ref{weak_betti} that $F_{\{1\}}$, $F_{x_4}$ and $F_{x_4^6}$ are the Markov fibers.
 
Let $S$ be a cardinality-minimal generating set of $I_L$. According to Theorem~\ref{main_thm}, $|S_{\overline F_{\{1\}}}|=\rank(L_{pure})=2$ and $S_{\overline F_{\{1\}}}$ must generate $I_{L_{pure}}=\langle 1-x_1^5,1-x_1x_2\rangle$. We note that by Theorem~\ref{Markov_pure} there are infinitely many choices for binomials $B_1, B_2$ that generate $I_{L_{pure}}$, and one way of constructing infinitely many such binomials is using Theorem~\ref{basis_primitive_lattice} with Corollary~\ref{cor_basis_pure}. We remark that the multiset $\{F_{B_1}, F_{B_2}\}$ equals $\{ F_{\{1\}}, F\}$, where $F$ can be any of the five fibers of ${\overline F_{\{1\}}}$.

Consider now the Markov fiber $F_{x_4}$. It is an easy exercise that 
\[
G(M_{F_{x_4}})=\{x_2^2x_3,x_1^3x_3,x_4\}.
\]
Thus $G(M_{F_{x_4}})^{\sigma}=\{x_4,x_3\}$. Since $I_{L,{<\overline F_{x_4}}}=I_{L_{pure}}$ it is immediate that $G_{\overline F_{x_4}}$ consists of two isolated vertices. Thus $t(\overline F_{x_4})=2$ and exactly one binomial is needed to construct a spanning tree of $\Gamma_{\overline F_{x_4}}$. To obtain a cardinality-minimal generating set, according to Theorem~\ref{main_thm}, we need to add to $S_{\overline F_{\{1\}}}$ a binomial $B_3=\pm (x^{\ub}-x^{\vb})$ such that $x^{\ub}\in x_3F_{x_1^{3+i}}$ and  $x^{\vb}\in x_4F_{x_1^i}$, where $0\leq i\leq 4$. For example, any of the binomials $x_1^2x_4 - x_3$ and $x_2x_3 - x_1^{2016}x_4$ can be chosen as $B_3$.

Next consider the Markov fiber $F'=F_{x_4^6}$, for which $\overline{F'}=\{ F_{x_4^6}, F_{x_1x_4^6}, F_{x_5}\}$. Based on Lemma~\ref{vertices_equiv_fibers_lem} and as in Example~\ref{graph_G_F} we have
\[
G(M_{F'})^{\sigma}=\{x_5,x_3^6,x_3^5x_4,x_3^4x_4^2,x_3^3x_4^3,x_3^2x_4^4,x_3x_4^5,x_4^6\}.
\]
Since the graph $G_{\overline F'}$ consists of two connected components (see Example~\ref{graph_G_F}), we need  one more binomial $B_4=\pm(x^{\wb_1}-x^{\wb_2})$ so that $x^{\wb_1}\in x_5F_{x_1^{i}}$ and $x^{\wb_2}$  is in the union of the sets  $x_3^6F_{x_1^{4+i}}$, $x_3^5x_4F_{x_1^{1+i}}$, $x_3^4x_4^2F_{x_1^{3+i}}$, $x_3^3x_4^3F_{x_1^{i}}$, $x_3^2x_4^4F_{x_1^{2+i}}$, $x_3x_4^5F_{x_1^{4+i}}$, $x_4^6F_{x_1^{1+i}}$,   $0\leq i\leq 4$. For example any of the binomials $x_3^2x_4^4 - x_2^2x_5$ and $x_1^2x_2^{2015}x_3x_4^5 - x_1^3x_5$ can be chosen as $B_4$.

According to Theorem~\ref{main_thm}, $S=\{B_1, B_2, B_3, B_4\}$ is a cardinality-minimal generating set of $I_L$. Thus $\mu(I_L)=4$ and $I_L$ is a binomial complete intersection. This also follows immediately from Corollary~\ref{char_binom_compl_inter_cor}, since $L^\sigma$ is generated by $v_1^\sigma$, $v_2^\sigma$ and  the matrix 
\[
 \begin{bmatrix} 1&-1&0\cr 6&0&-1\end{bmatrix}
\]
is mixed dominating.


\medskip

{\bf Acknowledgment}. The authors would like to thank Ezra Miller for useful discussions on binomial fibers, back in 2009, that initiated this work. Also, we thank the anonymous referee for several helpful suggestions on the terminology of this paper and for simplifying the proofs of Propositions~\ref{prim_dir}, \ref{lemma} and Theorems~\ref{primitive_bases}, \ref{basis_primitive_lattice}. This paper was partially written during the visit of the second and third author at the University of Ioannina. The third author was supported by a Romanian grant awarded by UEFISCDI, project number $83/2010$, PNII-RU code TE$\_46/2010$, program Human Resources, ``Algebraic modeling of some combinatorial objects and computational applications''.

\end{document}